\documentclass[11pt]{article}

\usepackage[T1, T2A]{fontenc}
\usepackage[english]{babel}
\usepackage[mathletters]{ucs} 
\usepackage[utf8x]{inputenc} 

\usepackage{adjustbox,amsbsy,amsfonts,amsmath,amssymb,amsthm,anyfontsize,authblk,bbding,bbm,bm,booktabs,braket,enumerate,eurosym,fancyvrb,float,graphicx,geometry,grffile,longtable,mathrsfs,multirow,textcomp,upquote,xcolor,tikz,tikz-cd,ytableau}

\usepackage{array,bigints,changepage,cite,comment,dsfont,egothic,enumerate,exscale,filecontents,latexsym,mathtools,relsize,rotating,stmaryrd}

\usepackage{titlesec}
\titlelabel{(\thetitle)\quad}

\newtheoremstyle{mystyle}
{}
{}
{\itshape}
{}
{\bfseries}
{.}
{.5em}
{\thmname{#1}\thmnumber{ #2}\thmnote{ (#3)}}

\theoremstyle{mystyle}
\newtheorem{definition}{Definition}
\newtheorem*{definition*}{Definition}

\newtheorem*{proposition*}{Proposition}
\newtheorem{lemma}{Lemma}
\newtheorem*{lemma*}{Lemma}
\newtheorem{theorem}{Theorem}
\newtheorem*{theorem*}{Theorem}
\newtheorem{corollary}{Corollary}
\newtheorem*{corollary*}{Corollary}

\newcommand{\sixptr}{{\scalebox{.5}{\text{\SixFlowerPetalRemoved}}}}

\newcommand*\circled[1]{\tikz[baseline=(char.base)]{
\node[shape=circle, draw, inner sep=1pt] (char) {#1};}}

\usepackage{nicefrac}

\usetikzlibrary{shapes,fit,matrix,chains,positioning,decorations.pathreplacing,arrows,positioning,circuits,decorations.markings}
\usetikzlibrary{calc}

\makeatletter

\usepackage[font=small,labelfont=bf, figurename=Fig.]{caption} 

\usepackage[font=small,labelfont=bf, tablename=Tab.]{caption}

\AtBeginDocument{%
} 
 
\usepackage[inline]{enumitem} 
\usepackage[normalem]{ulem} 

\definecolor{ansi-black}{HTML}{3E424D}
\definecolor{ansi-black-intense}{HTML}{282C36}
\definecolor{ansi-red}{HTML}{E75C58}
\definecolor{ansi-red-intense}{HTML}{B22B31}
\definecolor{ansi-green}{HTML}{00A250}
\definecolor{ansi-green-intense}{HTML}{007427}
\definecolor{ansi-yellow}{HTML}{DDB62B}
\definecolor{ansi-yellow-intense}{HTML}{B27D12}
\definecolor{ansi-blue}{HTML}{208FFB}
\definecolor{ansi-blue-intense}{HTML}{0065CA}
\definecolor{ansi-magenta}{HTML}{D160C4}
\definecolor{ansi-magenta-intense}{HTML}{A03196}
\definecolor{ansi-cyan}{HTML}{60C6C8}
\definecolor{ansi-cyan-intense}{HTML}{258F8F}
\definecolor{ansi-white}{HTML}{C5C1B4}
\definecolor{ansi-white-intense}{HTML}{A1A6B2}
\definecolor{ansi-hefault-inverse-fg}{HTML}{FFFFFF}
\definecolor{ansi-hefault-inverse-bg}{HTML}{000000}

\DefineVerbatimEnvironment{Highlighting}{Verbatim}{commandchars=\\\{\}}

\let\Oldtex\TeX
\let\Oldlatex\LaTeX
\renewcommand{\TeX}{\textrm{\Oldtex}}
\renewcommand{\LaTeX}{\textrm{\Oldlatex}}

\makeatletter
\def\PY@reset{\let\PY@it=\relax \let\PY@bf=\relax%
\let\PY@ul=\relax \let\PY@tc=\relax%
\let\PY@bc=\relax \let\PY@ff=\relax}
\def\PY@tok#1{\csname PY@tok@#1\endcsname}
\def\PY@toks#1+{\ifx\relax#1\empty\else%
\PY@tok{#1}\expandafter\PY@toks\fi}
\def\PY@do#1{\PY@bc{\PY@tc{\PY@ul{%
\PY@it{\PY@bf{\PY@ff{#1}}}}}}}
\def\PY#1#2{\PY@reset\PY@toks#1+\relax+\PY@do{#2}}

\expandafter\def\csname PY@tok@w\endcsname{\def\PY@tc##1{\textcolor[rgb]{0.73,0.73,0.73}{##1}}}
\expandafter\def\csname PY@tok@c\endcsname{\let\PY@it=\textit\def\PY@tc##1{\textcolor[rgb]{0.25,0.50,0.50}{##1}}}
\expandafter\def\csname PY@tok@cp\endcsname{\def\PY@tc##1{\textcolor[rgb]{0.74,0.48,0.00}{##1}}}
\expandafter\def\csname PY@tok@k\endcsname{\let\PY@bf=\textbf\def\PY@tc##1{\textcolor[rgb]{0.00,0.50,0.00}{##1}}}
\expandafter\def\csname PY@tok@kp\endcsname{\def\PY@tc##1{\textcolor[rgb]{0.00,0.50,0.00}{##1}}}
\expandafter\def\csname PY@tok@kt\endcsname{\def\PY@tc##1{\textcolor[rgb]{0.69,0.00,0.25}{##1}}}
\expandafter\def\csname PY@tok@o\endcsname{\def\PY@tc##1{\textcolor[rgb]{0.40,0.40,0.40}{##1}}}
\expandafter\def\csname PY@tok@ow\endcsname{\let\PY@bf=\textbf\def\PY@tc##1{\textcolor[rgb]{0.67,0.13,1.00}{##1}}}
\expandafter\def\csname PY@tok@nb\endcsname{\def\PY@tc##1{\textcolor[rgb]{0.00,0.50,0.00}{##1}}}
\expandafter\def\csname PY@tok@nf\endcsname{\def\PY@tc##1{\textcolor[rgb]{0.00,0.00,1.00}{##1}}}
\expandafter\def\csname PY@tok@nc\endcsname{\let\PY@bf=\textbf\def\PY@tc##1{\textcolor[rgb]{0.00,0.00,1.00}{##1}}}
\expandafter\def\csname PY@tok@nn\endcsname{\let\PY@bf=\textbf\def\PY@tc##1{\textcolor[rgb]{0.00,0.00,1.00}{##1}}}
\expandafter\def\csname PY@tok@ne\endcsname{\let\PY@bf=\textbf\def\PY@tc##1{\textcolor[rgb]{0.82,0.25,0.23}{##1}}}
\expandafter\def\csname PY@tok@nv\endcsname{\def\PY@tc##1{\textcolor[rgb]{0.10,0.09,0.49}{##1}}}
\expandafter\def\csname PY@tok@no\endcsname{\def\PY@tc##1{\textcolor[rgb]{0.53,0.00,0.00}{##1}}}
\expandafter\def\csname PY@tok@nl\endcsname{\def\PY@tc##1{\textcolor[rgb]{0.63,0.63,0.00}{##1}}}
\expandafter\def\csname PY@tok@ni\endcsname{\let\PY@bf=\textbf\def\PY@tc##1{\textcolor[rgb]{0.60,0.60,0.60}{##1}}}
\expandafter\def\csname PY@tok@na\endcsname{\def\PY@tc##1{\textcolor[rgb]{0.49,0.56,0.16}{##1}}}
\expandafter\def\csname PY@tok@nt\endcsname{\let\PY@bf=\textbf\def\PY@tc##1{\textcolor[rgb]{0.00,0.50,0.00}{##1}}}
\expandafter\def\csname PY@tok@nd\endcsname{\def\PY@tc##1{\textcolor[rgb]{0.67,0.13,1.00}{##1}}}
\expandafter\def\csname PY@tok@s\endcsname{\def\PY@tc##1{\textcolor[rgb]{0.73,0.13,0.13}{##1}}}
\expandafter\def\csname PY@tok@sd\endcsname{\let\PY@it=\textit\def\PY@tc##1{\textcolor[rgb]{0.73,0.13,0.13}{##1}}}
\expandafter\def\csname PY@tok@si\endcsname{\let\PY@bf=\textbf\def\PY@tc##1{\textcolor[rgb]{0.73,0.40,0.53}{##1}}}
\expandafter\def\csname PY@tok@se\endcsname{\let\PY@bf=\textbf\def\PY@tc##1{\textcolor[rgb]{0.73,0.40,0.13}{##1}}}
\expandafter\def\csname PY@tok@sr\endcsname{\def\PY@tc##1{\textcolor[rgb]{0.73,0.40,0.53}{##1}}}
\expandafter\def\csname PY@tok@ss\endcsname{\def\PY@tc##1{\textcolor[rgb]{0.10,0.09,0.49}{##1}}}
\expandafter\def\csname PY@tok@sx\endcsname{\def\PY@tc##1{\textcolor[rgb]{0.00,0.50,0.00}{##1}}}
\expandafter\def\csname PY@tok@m\endcsname{\def\PY@tc##1{\textcolor[rgb]{0.40,0.40,0.40}{##1}}}
\expandafter\def\csname PY@tok@gh\endcsname{\let\PY@bf=\textbf\def\PY@tc##1{\textcolor[rgb]{0.00,0.00,0.50}{##1}}}
\expandafter\def\csname PY@tok@gu\endcsname{\let\PY@bf=\textbf\def\PY@tc##1{\textcolor[rgb]{0.50,0.00,0.50}{##1}}}
\expandafter\def\csname PY@tok@gd\endcsname{\def\PY@tc##1{\textcolor[rgb]{0.63,0.00,0.00}{##1}}}
\expandafter\def\csname PY@tok@gi\endcsname{\def\PY@tc##1{\textcolor[rgb]{0.00,0.63,0.00}{##1}}}
\expandafter\def\csname PY@tok@gr\endcsname{\def\PY@tc##1{\textcolor[rgb]{1.00,0.00,0.00}{##1}}}
\expandafter\def\csname PY@tok@ge\endcsname{\let\PY@it=\textit}
\expandafter\def\csname PY@tok@gs\endcsname{\let\PY@bf=\textbf}
\expandafter\def\csname PY@tok@gp\endcsname{\let\PY@bf=\textbf\def\PY@tc##1{\textcolor[rgb]{0.00,0.00,0.50}{##1}}}
\expandafter\def\csname PY@tok@go\endcsname{\def\PY@tc##1{\textcolor[rgb]{0.53,0.53,0.53}{##1}}}
\expandafter\def\csname PY@tok@gt\endcsname{\def\PY@tc##1{\textcolor[rgb]{0.00,0.27,0.87}{##1}}}
\expandafter\def\csname PY@tok@err\endcsname{\def\PY@bc##1{\setlength{\fboxsep}{0pt}\fcolorbox[rgb]{1.00,0.00,0.00}{1,1,1}{\strut ##1}}}
\expandafter\def\csname PY@tok@kc\endcsname{\let\PY@bf=\textbf\def\PY@tc##1{\textcolor[rgb]{0.00,0.50,0.00}{##1}}}
\expandafter\def\csname PY@tok@kd\endcsname{\let\PY@bf=\textbf\def\PY@tc##1{\textcolor[rgb]{0.00,0.50,0.00}{##1}}}
\expandafter\def\csname PY@tok@kn\endcsname{\let\PY@bf=\textbf\def\PY@tc##1{\textcolor[rgb]{0.00,0.50,0.00}{##1}}}
\expandafter\def\csname PY@tok@kr\endcsname{\let\PY@bf=\textbf\def\PY@tc##1{\textcolor[rgb]{0.00,0.50,0.00}{##1}}}
\expandafter\def\csname PY@tok@bp\endcsname{\def\PY@tc##1{\textcolor[rgb]{0.00,0.50,0.00}{##1}}}
\expandafter\def\csname PY@tok@fm\endcsname{\def\PY@tc##1{\textcolor[rgb]{0.00,0.00,1.00}{##1}}}
\expandafter\def\csname PY@tok@vc\endcsname{\def\PY@tc##1{\textcolor[rgb]{0.10,0.09,0.49}{##1}}}
\expandafter\def\csname PY@tok@vg\endcsname{\def\PY@tc##1{\textcolor[rgb]{0.10,0.09,0.49}{##1}}}
\expandafter\def\csname PY@tok@vi\endcsname{\def\PY@tc##1{\textcolor[rgb]{0.10,0.09,0.49}{##1}}}
\expandafter\def\csname PY@tok@vm\endcsname{\def\PY@tc##1{\textcolor[rgb]{0.10,0.09,0.49}{##1}}}
\expandafter\def\csname PY@tok@sa\endcsname{\def\PY@tc##1{\textcolor[rgb]{0.73,0.13,0.13}{##1}}}
\expandafter\def\csname PY@tok@sb\endcsname{\def\PY@tc##1{\textcolor[rgb]{0.73,0.13,0.13}{##1}}}
\expandafter\def\csname PY@tok@sc\endcsname{\def\PY@tc##1{\textcolor[rgb]{0.73,0.13,0.13}{##1}}}
\expandafter\def\csname PY@tok@dl\endcsname{\def\PY@tc##1{\textcolor[rgb]{0.73,0.13,0.13}{##1}}}
\expandafter\def\csname PY@tok@s2\endcsname{\def\PY@tc##1{\textcolor[rgb]{0.73,0.13,0.13}{##1}}}
\expandafter\def\csname PY@tok@sh\endcsname{\def\PY@tc##1{\textcolor[rgb]{0.73,0.13,0.13}{##1}}}
\expandafter\def\csname PY@tok@s1\endcsname{\def\PY@tc##1{\textcolor[rgb]{0.73,0.13,0.13}{##1}}}
\expandafter\def\csname PY@tok@mb\endcsname{\def\PY@tc##1{\textcolor[rgb]{0.40,0.40,0.40}{##1}}}
\expandafter\def\csname PY@tok@mf\endcsname{\def\PY@tc##1{\textcolor[rgb]{0.40,0.40,0.40}{##1}}}
\expandafter\def\csname PY@tok@mh\endcsname{\def\PY@tc##1{\textcolor[rgb]{0.40,0.40,0.40}{##1}}}
\expandafter\def\csname PY@tok@mi\endcsname{\def\PY@tc##1{\textcolor[rgb]{0.40,0.40,0.40}{##1}}}
\expandafter\def\csname PY@tok@il\endcsname{\def\PY@tc##1{\textcolor[rgb]{0.40,0.40,0.40}{##1}}}
\expandafter\def\csname PY@tok@mo\endcsname{\def\PY@tc##1{\textcolor[rgb]{0.40,0.40,0.40}{##1}}}
\expandafter\def\csname PY@tok@ch\endcsname{\let\PY@it=\textit\def\PY@tc##1{\textcolor[rgb]{0.25,0.50,0.50}{##1}}}
\expandafter\def\csname PY@tok@cm\endcsname{\let\PY@it=\textit\def\PY@tc##1{\textcolor[rgb]{0.25,0.50,0.50}{##1}}}
\expandafter\def\csname PY@tok@cpf\endcsname{\let\PY@it=\textit\def\PY@tc##1{\textcolor[rgb]{0.25,0.50,0.50}{##1}}}
\expandafter\def\csname PY@tok@c1\endcsname{\let\PY@it=\textit\def\PY@tc##1{\textcolor[rgb]{0.25,0.50,0.50}{##1}}}
\expandafter\def\csname PY@tok@cs\endcsname{\let\PY@it=\textit\def\PY@tc##1{\textcolor[rgb]{0.25,0.50,0.50}{##1}}}

\makeatother

\definecolor{incolor}{rgb}{0.0, 0.0, 0.5}
\definecolor{outcolor}{rgb}{0.545, 0.0, 0.0}

\usepackage[none]{hyphenat}
\emergencystretch=\maxdimen
\usepackage{scalerel,stackengine}
\def\scaleone#1{\vcenter{\hbox{\scaleto[3ex]{\displaystyle\mathds{1}}{#1}}}}

\usepackage{silence}
\usepackage[square,sort&compress,comma,numbers]{natbib}

\makeatletter
\renewcommand{\@biblabel}[1]{[#1]\hfill}

\title{\vskip-1.15in\fontsize{14}{19}\selectfont 
Asymptotic Shape of Quantum Markov Semigroups for Compact Uniform Trees\vskip-.125in
} 

\author{Margarita Belova$^{\,(1\dagger)}$ and Matthew Bernard$^{\,(2\dagger)}$\vskip-.0in
	
{\small $^{(1\dagger)\,}$margarita.s.belova@phystech.edu, 
\,$^{(2\dagger)\,}$mattb@berkeley.edu}} 

\affil{\vskip-.1in 
{\small $^{(1\dagger)\,}$Moscow Inst$.$ of Physics and Technology 
(MIPT)$,$ Dolgoprudny$,$ Russia\vskip.0in
$^{(2\dagger)\,}$Advanced Computational Biology Center (ACBC)$,$ Berkeley$,$ CA$,$ USA\vskip-.55in $\,$
}}

\date{
}

\begin{document}

\maketitle

\vskip-.45in$\,$

\begin{abstract}\vskip-.05in
\noindent We give locally finite Markov trees in $L^p$-compact$,$ separable Hilbert$,$ supersymmetric process$:$ $[0,\infty)\!\times\!\mathbb{R}^{\lvert\mathcal{A}^{\otimes m}\rvert}/\mathcal{A}^{\otimes m}$ on quantum ${\rm U}(\lvert\mathcal{A}^{\otimes m}\rvert)$ semigroups$.$ In full automorphism group ${\rm Aut}({\rm\bf T})$ of modular subgroup$,$ asymptotic-ergodicity is entropy-worthy $\mathbb{R}$ shape for uniform partition$.$ 
\\[.05in]
\end{abstract}\vskip-.3in $\,$

\subsection*{$\bm{GL(n)}$ theory for discrete-, continuous-time, quantum $\bm{{\rm U}(\lvert\mathcal{A}^{\otimes m}\rvert)}$ Markov trees} $\,$\vskip-.2in

\begin{figure}[H]
\centering
\begin{tikzpicture}[scale=4.0]
\draw (0,0) node[anchor = east] {$\Phi_0$} 
-- (0.8,0) node[anchor = west] {$\Phi_1$};
\draw (0.8,0) -- (0.9,0.55) node [anchor = west] {$\Phi_2$};
\draw (0.8,0) -- (0.4,.69282) node [anchor = south] {$\Phi_3$};
\draw (0,0) -- (0.4,.69282);
\draw (0.9,0.55) -- (0.4,.69282);
\draw[dashed]  (0,0)--(0.9,0.55);
\end{tikzpicture}\vskip-.1in
\caption{$(<\infty)$-mixture space of trees$:$ $4$-mixture space$,$ namely$,$ tetrahedron ($3$-simplex)$.$}
\label{simplices}\vskip-.1in
\end{figure}
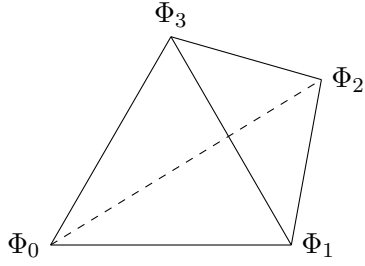

\noindent Adhering $(\cdot,\cdot)_\Phi,$ $\mathcal{X}\!\equiv\!L^p(\mathcal{X},\Phi)$-compact space $(\Omega,\mathcal{F},{\rm\bf P})$ of usual $[0,\infty)\!\times\!\mathbb{R}^{\lvert\mathcal{A}^{\otimes m}\rvert}$ \!Markov trees$,$ is $\mathcal{X}$-acting process$,$ positive semi-definite $H$ set $\mathcal{P}$ of density (normalized)$,$ subnormalized states$:$
\\[-.25in]
\begin{align} \hskip-.1in
\{{\rm\bf T}\!:\,{tr}({\rm\bf T})\!=1\} 
\,\subseteq\,\mathcal{P}, 
~ \text{ resp$.$ } ~ 
\{{\rm\bf T}\!:\,{tr}({\rm\bf T})\!\leqslant\!1\} 
\,\subseteq\,\mathcal{P} 
\end{align} 
\\[-.275in]
in supersymmetric basis $\mathcal{A}^{\otimes m}$ representation of the real Lie group$,$ unitary group ${\rm U}(\lvert\mathcal{A}^{\otimes m}\rvert)$ for $\nu$-partite Hilbert space $\mathcal{X}_{X_0\cdots X_{\nu-1}} 
\!= \mathcal{X}_0\!\oplus\!\cdots\!\oplus\!\mathcal{X}_{\nu-1},$ and affine functional ${\rm\bf\Gamma}\!\!\longrightarrow\!\Phi_{\rm\bf\Gamma}(F),$ $F\!\in\!\mathcal{F}\!:$ 
\\[-.225in]
\begin{align}\hskip-.075in 
{\rm\bf\Gamma}_{\!X_\xi} 
\!={\rm tr}_{X_\xi}\!\big({\rm\bf\Gamma}_{\!X_0\cdots X_{\nu-1}}\big), 
\,\forall\,\xi\!=\!0,\ldots,\nu\!-\!1 
\hskip.025in \Big| \hskip.025in 
\Phi_{{\rm\bf\Gamma}_{\!x}}(F)\!\neq\!\Phi_{{\rm\bf\Gamma}_{\!y}}(F);
\text{ i.e$.$ } 
\exists\,{\rm\bf\Gamma}\!\!\longrightarrow\!\Phi_{{\rm\bf\Gamma}},
\,\Phi_{{\rm\bf\Gamma}_{\!x}}\!\not\equiv\!\Phi_{{\rm\bf\Gamma}_{\!y}},
\,\forall\,{\rm\bf\Gamma}_{\!x}\!\neq\!{\rm\bf\Gamma}_{\!y} 
\end{align}
\\[-.25in]
where ${\rm\bf\Gamma}_{\!X_0},\ldots,{\rm\bf\Gamma}_{\!X_{\nu-1}}$ are reduced states of unique mixed state ${\rm\bf\Gamma}_{\!X_0\cdots X_{\nu-1}}\in\mathcal{P}_{\!X_0\cdots X_{\nu-1}}$ on compact (closed and bounded) space$.$ For events $\mathbb{X}\!\ni\!x\!\in\!\mathcal{X},$ measurable subsets $\mathbb{X}\!\subset\!\mathcal{X}$ all form $\sigma$-field $\mathcal{F}(\mathcal{X})$ of measurable space $(\mathcal{X},\mathcal{F})$ for all outcomes$.$ In real$,$ Euclidean $\lvert\mathcal{A}^{\otimes m}\rvert$-dimensional space $\mathbb{R}^{\lvert\mathcal{A}^{\otimes m}\rvert},$ $\mathcal{X}$ is a domain with Borel $\sigma$-field generated by multi-dimensional open sets (intervals)$.$ Main result here is limiting affine bilinear mapping $\gamma\!:\mathcal{P}\!\longrightarrow\!\Phi;$ by $\Phi\!:\mathbb{R}^{\lvert\mathcal{A}^{\otimes m}\rvert}\!\longrightarrow\!\mathbb{R}_+,$ $\forall\,\left\lVert{\rm\bf T}\!-\!{\rm\bf T}^{\sixptr}\right\rVert_{\mathcal{L}(L^{p\geqslant1}(\gamma))}\!\leqslant\!\varepsilon\!:$
\\[-.225in]
\begin{align}\label{mu}\hskip-.075in 
\Phi_{\!X}= 
\hskip-.35in \sum_{x\,\in\,(\mathcal{A}^{\otimes m}\,\cup\,\{{\rm\bf I}_{\lvert\mathcal{A}^{\otimes m}\rvert}\})} 
\hskip-.45in {\rm\bf P}_{\!\!x}\,tr(x^\dagger x \,{\rm\bf\Gamma}_{\!x}), 
\hskip.025in \text{ resp$.$ }  
\hskip.025in \Phi_{\!X}\!=\!\!\bigintsss_{\mathbb{R}^{\lvert\mathcal{A}^{\otimes m}\rvert}} 
\hskip-.25in \det\!\bigg( 
\hskip-.00in \sum_{x\,\in\,(\mathcal{A}^{\otimes m}\,\cup\,\{{\rm\bf I}_{\lvert\mathcal{A}^{\otimes m}\rvert}\})} 
\hskip-.45in x^\dagger x 
\hskip.050in \exp\!\Big(\!\!-iH_{\rm\bf\Gamma}(\theta_x) 
\hskip-.00in \Big) 
\hskip-.025in \bigg) 
\hskip.025in {\rm\bf P}_{\!t}\big(\vec{\theta_x}\in d\vec{\theta_x}\big) 
\end{align} 
\\[-.475in]
\begin{align} 
\hskip-.0in x\equiv x_{\xi_0}\!\cdots x_{\xi_{m-1}} 
\!= \left\lvert x_{\xi_0} \right\rangle 
\!\otimes\cdots\otimes\!\left\lvert x_{\xi_{m-1}}\right\rangle; 
\hskip.15in x^\dagger x \equiv \ket{x_{i1}^\dagger}\!\bra{x_{1j}};
\hskip.15in \overline{x^{\rm T}}\!=x^\dagger; 
\end{align} 
\\[-.55in]
\begin{align}\hskip-.0in {\rm\bf T}_{\!X} \!= 
\hskip-.25in \sum_{x\,\in\,(\mathcal{A}^{\otimes m}\,\cup\,\{{\rm\bf I}_{\lvert\mathcal{A}^{\otimes m}\rvert}\})} 
\hskip-.45in {\rm\bf P}_{\!\!x} ~ x^\dagger x, 
\hskip.125in \text{ resp$.$ }  
\hskip.125in {\rm\bf T}_{\!X} \!=  
\hskip-.25in \sum_{x\,\in\,(\mathcal{A}^{\otimes m}\,\cup\,\{{\rm\bf I}_{\lvert\mathcal{A}^{\otimes m}\rvert}\})} 
\hskip-.45in x^\dagger x ~ \exp\!\Big(\!\!-\delta(x)\Big) 
\hskip.025in {\rm\bf P}_{\!t}\big(\vec{\theta_x}\in d\vec{\theta_x}\big). 
\end{align} 
\\[-.175in] 
Dirac-delta $\delta(x)$ thus assumes zero Hamiltonian of the generic Lagrangian guage fixings (quantum  gravity)$.$ Diagonal state ${\rm\bf\Gamma}$ is on eigenfunction as unique probabilities $\{{\rm\bf P}_{\!\!x}\!=\!{\rm\bf P}_{\!X}(x)\!=\!{\rm\bf P}(X\!\!=\!x)\}_{x \,\in\, \mathcal{X}},$ resp$.$ ${\rm\bf P}_{\!t}$ on i.i.d $(\delta_{ij})$ random $X$ from set $\mathcal{X}.$ For all infinite dimensional analogue$:$ $(n\!+\!1)$-mixture is convex hull of a set of vertices$:$ $n$-simplex extreme points $\{\Phi_{\xi}\big|^{2^m}_{\xi=0}\}$ (Fig$.$~\ref{simplices})$;$ a simplex of mixture sequence $\{{\rm\bf\Gamma}\}\!\subset\!\mathcal{P}$ gets modular algebraic closure \cite{jc79,rb92,cgob19,bgiko19} of stellations and spherical tilings in weights $\{{\rm\bf P}_{\!\!x}\}$ convex mixture$,$ by the one-to-one map of $\mathcal{P}$ into convex subset of linear space$.$ The bipartite conditional states $\{{\rm\bf\Gamma}_{\!Q\mid X=\,x}\}_{x \,\in\, \mathcal{X}}$ operator is feasible as$:$ 
\\[-.2in] 
\begin{align}
\hskip-.0in {\rm\bf\Gamma}_{\!XQ} = 
\hskip-.25in \sum_{x\,\in\,(\mathcal{A}^{\otimes m}\,\cup\,\{{\rm\bf I}_{\lvert\mathcal{A}^{\otimes m}\rvert}\})} 
\hskip-.45in {\rm\bf P}_{\!\!x} 
~ x^\dagger x\otimes{\rm\bf\Gamma}_{\!Q\mid X=\,x} 
\hskip.1in \Big| ~ 
{\rm\bf\Gamma}_{\!XQ} = 
\hskip-.25in \sum_{x\,\in\,(\mathcal{A}^{\otimes m}\,\cup\,\{{\rm\bf I}_{\lvert\mathcal{A}^{\otimes m}\rvert}\})} 
\hskip-.45in x^\dagger x 
\hskip.050in \exp\!\Big(\!\!-iH_{{\rm\bf\Gamma}_{\!Q\mid X=\,x}}(\theta_x) 
\hskip-.00in \Big) 
\hskip.025in {\rm\bf P}_{\!t}\big(\vec{\theta_x}\in d\vec{\theta_x}\big) 
\end{align} 
\\[-.15in] 
For the \cite{defkz00,oss16,oss19} real tree $(\mathcal{T}_{{\rm\bf T},{\rm\bf T}^{\sixptr}}\,,\,\varepsilon^\ast)$ semigroup $({\rm\bf P}_{\!t})_{t\geqslant0}\!:$ Take locally finite uniform random tree as the universal cover of finite (graph) random partitions, with full automorphism group ${\rm Aut}({\rm\bf\Gamma})$ as locally compact topological group$.$ By strongly continuous one-parameter semigroups$,$ on infinite sequence of tree states ${\rm\bf\Gamma}$ process distribution $\{{\rm\bf P}_{\!t}\},$ thus$,$ ($\mathcal{P}, \gamma$) defines structural (Markov) property in equivalence class partition $[\widetilde{\rm\bf\Gamma}]$ on conditional probability of $x\!\in\!\mathbb{X}$ under initial condition of $\widetilde{\rm\bf\Gamma},$ explicitly by virtue of automorphic functions with respect to subgroup of modular group$,$ where the $(\cdot,\cdot)_\Phi$ affine bilinear map $\gamma$ acts as measure transforming states into distributions on $\mathcal{X},$ for all $\Phi.$
$\,$\vskip-.2in 

\begin{figure}[H]
\centering ~ \hskip-0.in 
\includegraphics[width=.15\textwidth,angle=0,origin=c]
{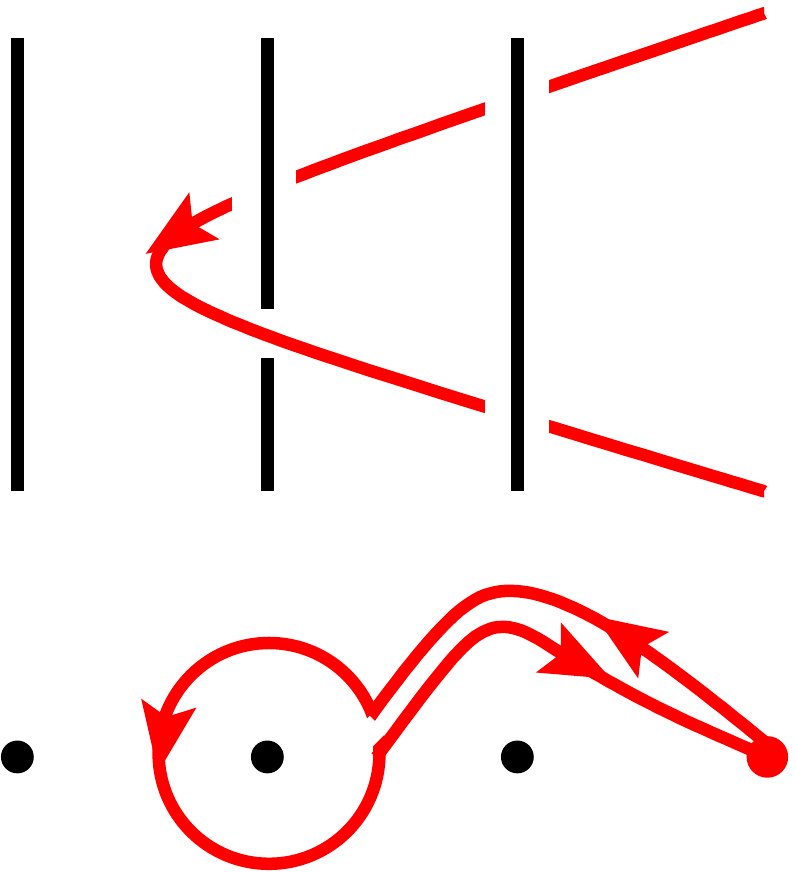}\hskip.1in 
\includegraphics[width=.45\textwidth,angle=0,origin=c]
{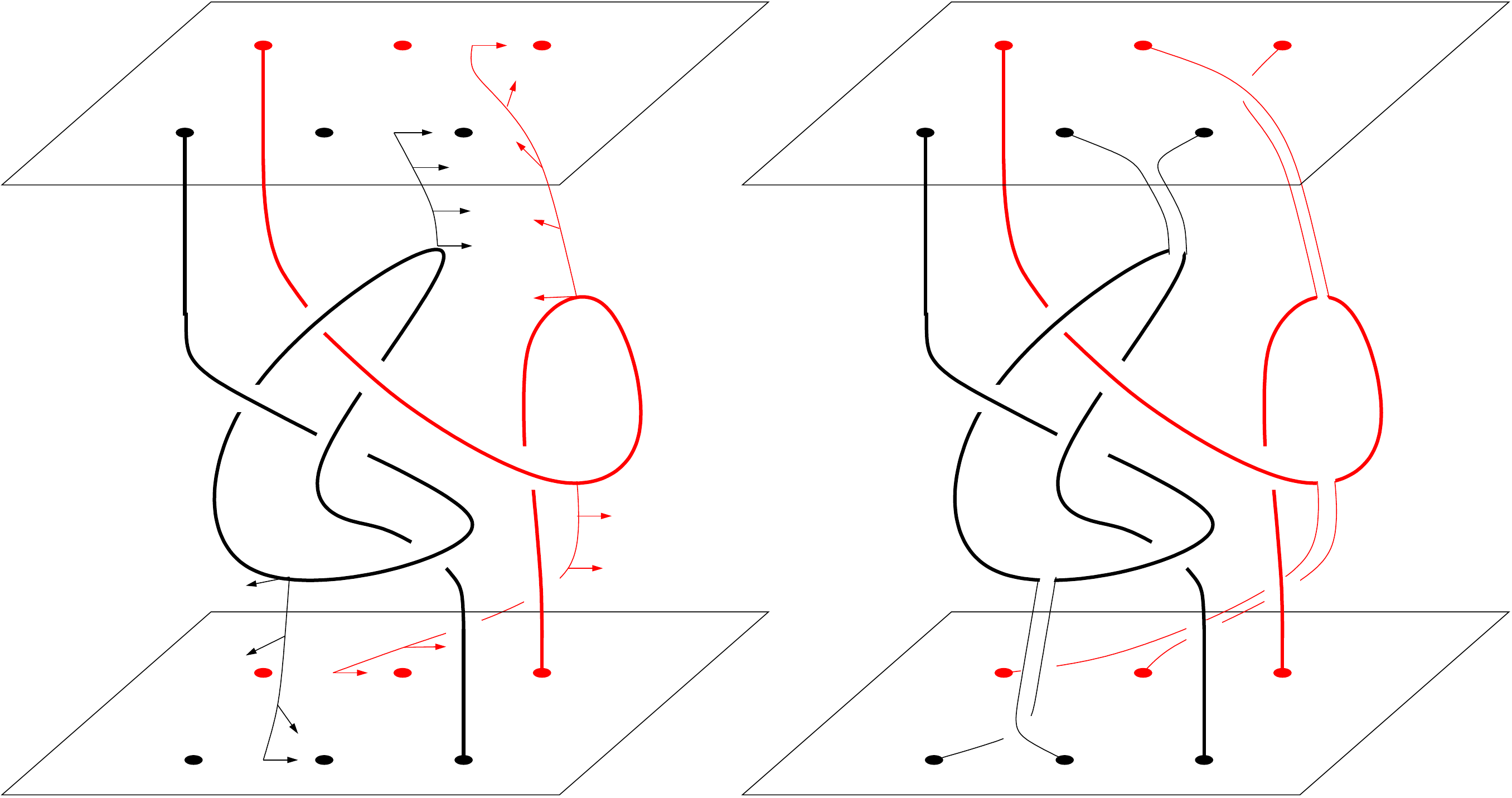}\hskip.1in 
\includegraphics[width=.25\textwidth,angle=0,origin=c]
{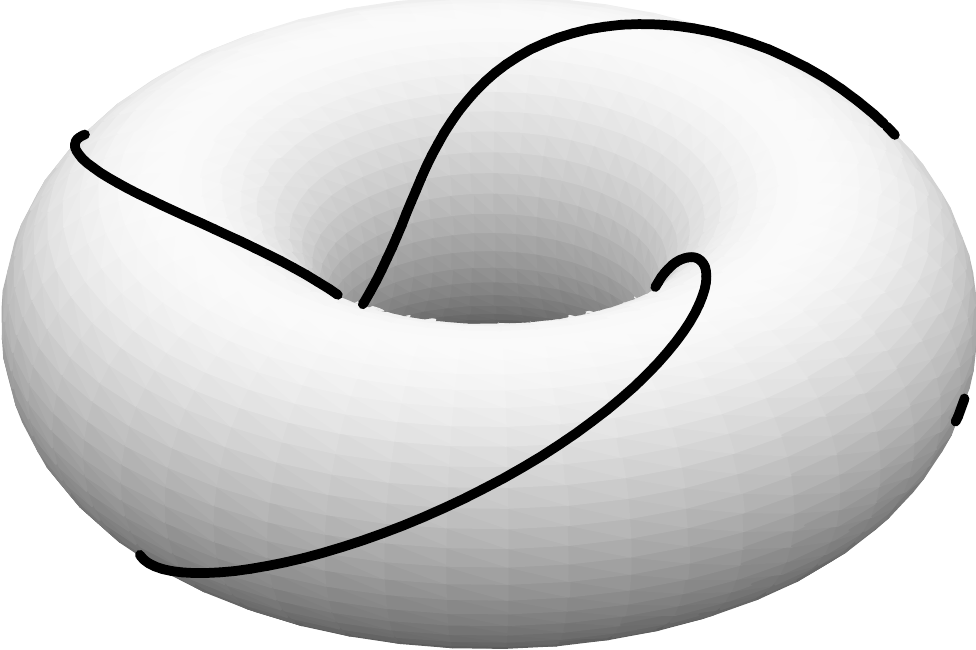}
\hskip-0in\,\\[-.05in]
\caption{Tree-isometry placed braid as string-link wound in composition of generators (permutations)$.$}\label{braid-string-link} \vskip-.075in 
\end{figure}

\begin{figure}[H]
\centering
\includegraphics[width=.65\textwidth]{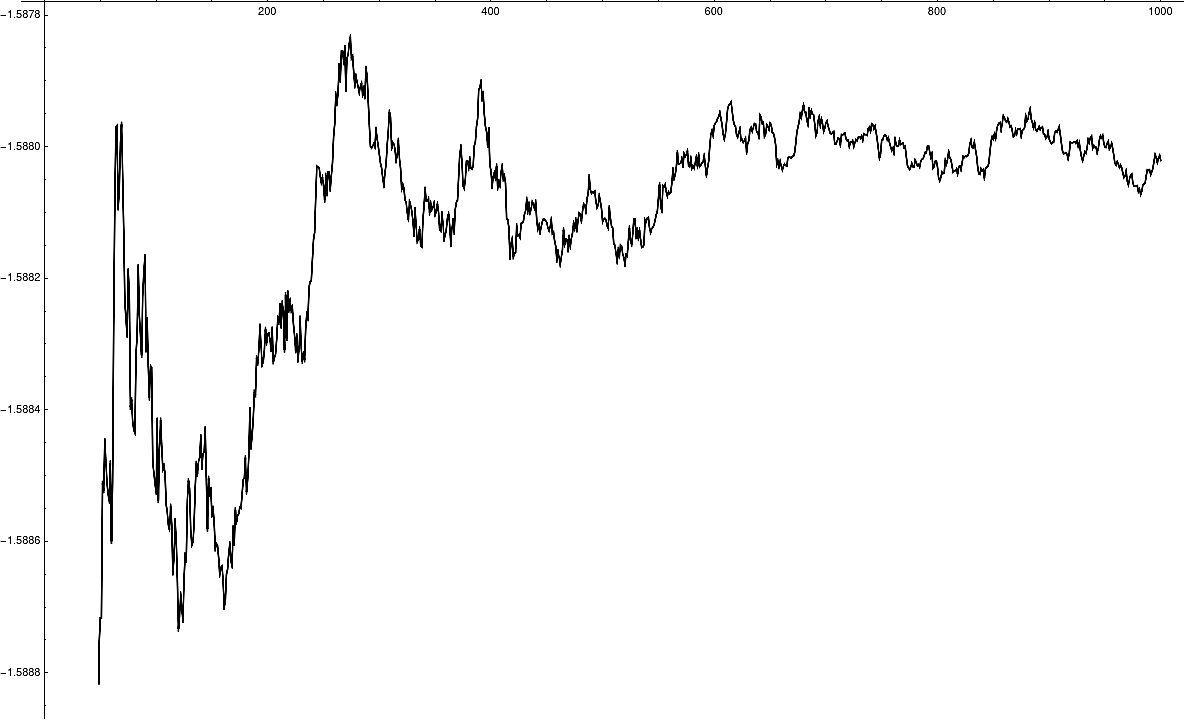}\\[-.05in]
\caption{minus entropy$,$ limiting asymptotic tree shape of no cut-off for affine density $f \!= 5 + \mathcal{N}(0,1).$}\label{2dentropy}  
\end{figure}

\begin{definition}\label{def1} 
Let $({\rm\bf P}_{\!\!x},x\!\in\!\mathcal{X})$ be a family of $\gamma$ probability measures for Markov process $(X_t)$ taking values in $\mathcal{X}$ set up on limit shape probability space $(\Omega,\mathcal{F});$ then ${\rm\bf P}_{\!\!x}$ is the $L^p(\gamma),$ $p\!\geqslant\!1,$ distribution law of random variable $X$ under initial condition $X_0\!=\!x.$ For $\sigma$-algebra $\mathbb{X}$ on state space $\mathcal{X}$ of $X\!\!:$ 
\vskip.075in
\begin{enumerate}[label=(\roman*),itemsep=4pt,parsep=-0pt,topsep=0pt,partopsep=0pt,labelsep=.2cm,labelwidth=.0cm,align=right,itemindent=-0.25cm]
\item $(t,\omega)\!\longmapsto\!\!X_t(\omega)$ is  $\mathbb{B}\otimes\mathcal{F}\!/\mathbb{X}$-measurable mapping of $[0,\infty)\!\times\!\Omega$ into $\mathcal{X}$\! by $[0,\infty)$ \!Borel $\sigma$-algebra $\mathbb{B}$
\item $x\longmapsto\!{\rm\bf P}_{\!\!x}\{F\}$ is $\mathbb{X}$-measurable for each $F\!\in\!\mathcal{F};$ in particular$,$ for topological space $(\mathcal{X},2^{[\mathcal{X}]}).$ 
\end{enumerate}
\end{definition}
$\,$\\[-.2in] 
{\it Remark.}~~Physical observation suggests $(\Omega,\mathcal{F})$ accepts a random variable $T_\alpha \!\mid\! \alpha\!>\!0,$ which$,$ under ${\rm\bf P}_{\!\!x} \mid x\in\mathcal{X},$ is independent of $X$ and has the mixing$,$ exponential distribution of parameter $\alpha.$ 
\\[-.2in]
\begin{proposition*}[Feynman path semigroup] 
For $f\!:\mathbb{R}^{\lvert\mathcal{A}^{\otimes m}\rvert}\!\longrightarrow\!\mathbb{R}_+/\mathcal{A}^{\otimes m}\!\cup\!\{{\rm\bf I}_{\lvert\mathcal{A}^{\otimes m}\rvert}\},$ set ${\rm\bf P}_{\!t}f\!\mid_{t\geqslant0}:$
\\[-.2in]
\begin{align}\hskip-.95in 
\underbrace{\int_{\mathbb{R}} \hskip-.05in\cdots\hskip-.05in 
\int_{\mathbb{R}}}_{\lvert\mathcal{A}^{\otimes m}\rvert\text{ times}} 
\hskip-.025in\frac{f(x_1,\ldots,x_{2^m})}{\sqrt{\det(2\pi t\Sigma)}} 
\exp\!\Big( 
\!\!-\!\frac{1}{2}\!\sum^{\lvert\mathcal{A}^{\otimes m}\rvert}_{i,j\,=\,1} 
\!(x_i\!-\!\mu_i)t^{-1}\Sigma_{ij}^{-1}(x_j\!-\!\mu_j)\!\Big) dx_1\!\cdots dx_{\lvert\mathcal{A}^{\otimes m}\rvert} 
\hskip.05in \mathlarger{\mathlarger{\Bigg|}}\hskip-.05in 
\begin{array}{l}
(\Sigma^{-1}_{ij}) = \Sigma^{-1}; ~ (\Sigma_{ij})\!=\!\Sigma
\\[.025in]
\Sigma_{i,j\,=\,i}=\Sigma_{ii} 
\\[.025in]
\Sigma_{i,j\,\neq\,i}=\rho_{ij}\sqrt{\Sigma_{ii}\Sigma_{jj}}\mid j\!>\!i 
\end{array}\label{semigroup} \hskip-.4in  
\end{align}
\\[-.2in]
and put ${\rm\bf P}_{\!0}f=f.$ Then $({\rm\bf P}_{\!t})_{t\geqslant0}$ is a semigroup$;$ moreover$,$ on Feynman path ``amplitude$.$'' 
\end{proposition*}
$\,$\\[-.2in] 
{\it Proof.}~~${\rm\bf P}_{\!\!s}\!\circ\!{\rm\bf P}_{\!t}\!=\!{\rm\bf P}_{\!\!s+t};$ and$,$ clearly$,$ the integrand is of the form $e^{-(i/\hbar)S[x]}.$ In addition$,$ ${\rm\bf P}_{\!t}f(x)$ exists since ${\rm\bf P}_{\!t}f(x)=\mathbb{E}[f(x\!+\!\sqrt{t}Z)],$ for $Z\sim\mathcal{N}[0,1],$ on space of $(\lvert\mathcal{A}^{\otimes m}\rvert\!\times\!\lvert\mathcal{A}^{\otimes m}\rvert)$ symmetric invertible positive semi-definite $\Sigma\!=\!A^T\!A$ invariant
with respect to action of orthogonal group (i.e$.$ $U\Sigma\overset{\text{d}}{=}\Sigma$)$;$ 
and$,$ $\Sigma = P D_\Sigma P^{-1}$ 
for diagonal $D_\Sigma$ of eigenvalues 
and invertible $P$ of eigenvectors$.$\hfill $\square$
\\[.1in]
{\it Remark.}~~By (\ref{double})$,$ integral (\ref{semigroup}) makes sense for $f\!\in\!L^1,$ resp$.$ $f\!\in\!L^\infty,$ and defines a quantum$,$ linear contraction operator $L^1\!\longrightarrow\!L^1,$ resp$.$ $L^\infty\!\longrightarrow\!L^\infty.$ Consequently (Riesz-Thorin \cite{ri27,th48})$,$ ${\rm\bf P}_{\!t}$ can be thought of as quantum$,$ linear contraction operator $L^p\!\longrightarrow\!L^p$ restricted to upper half-plane$,$ for all $1\!\leqslant\!p\!\leqslant\!\infty.$ For a heat semigroup$,$ by $f\!:\mathbb{R}^{\lvert\mathcal{A}^{\otimes m}\rvert}\!\longrightarrow\!\mathbb{R},$ we have$,$ 
\\[-.25in]
\begin{align}
{\rm\bf P}_{\!t}f(\mu)=\!\int_{\mathbb{R}^{\lvert\mathcal{A}^{\otimes m}\rvert}} 
\hskip-.000in f(x)\frac{1}{(2\pi t)^{\frac{1}{2}\lvert\mathcal{A}^{\otimes m}\rvert}} 
\exp\!\bigg(\!\!-
\frac{\left\lVert x-\mu \right\rVert^2}{2t}\bigg) dx 
\hskip.2in \Big| \hskip.05in {\rm\bf P}_{\!0}f=f.
\end{align}
\\[-.4in]
\begin{lemma}[tree process] 
For stopping time $t,$ tree-valued Markov process $(X_t)_{t\geqslant0},$ there exists ${\rm\bf P}_{\!t}$ random process$,$ i.e$.$ random function on sample functions defined by $t\longmapsto W(t,\omega)$ over sample path$,$ as uniform limit of interpolated $n$-state continuous functions$,$ on independent increment$:$
\\[-.25in]
\begin{align} 
\hskip-.0in W(t_1)-W(0), ~ W(t_2)-W(t_1), ~ \ldots, ~ W(t_n)-W(t_{n-1}), 
\hskip.2in 
\forall ~ \lvert\mathcal{A}^{\otimes m}\rvert=n\in\mathbb{N}. 
\hskip-.0in
\end{align}
\\[-.25in] 
with respect to diagonal-eigenvalues and states$,$ given by $\left\lVert\lambda-\lambda^{\ast}\right\lVert$ and $\left\lVert{\rm\bf T}-{\rm\bf T}^{\sixptr}\right\rVert_{\mathcal{L}(L^p(\gamma))}.$
\end{lemma}
$\,$\\[-.2in] 
{\it Proof.}~~The idea is to construct in $n$-box$,$ a time $[0,1],$ $n$-dimensional (i.e$.$ $n$-tuple state) process starting at $(0,\ldots,0)$ and interpolated on time intervals $1\!\leqslant\!\xi\!\leqslant\!f_\eta,$ $\forall\,\eta\!\geqslant\!0,$ as $n$-tuple random element in space $C[0,1]$ of continuous functions$,$ from ensemble of $n$ paths of respective time $[0,1]$ planar continuous-state process with state-space $\mathcal{X}$ starting at zero i.e$.$ from $n$-path ensemble of time $[0,1]$ planar continuous-state random walk of state-space $\mathcal{X}$ starting at zero$.$ \hfill $\square$ 
\\[-.65in]
\begin{figure}[H]
\centering ~ \hskip-0.in 
\includegraphics[width=.45\textwidth,angle=0,origin=c]
{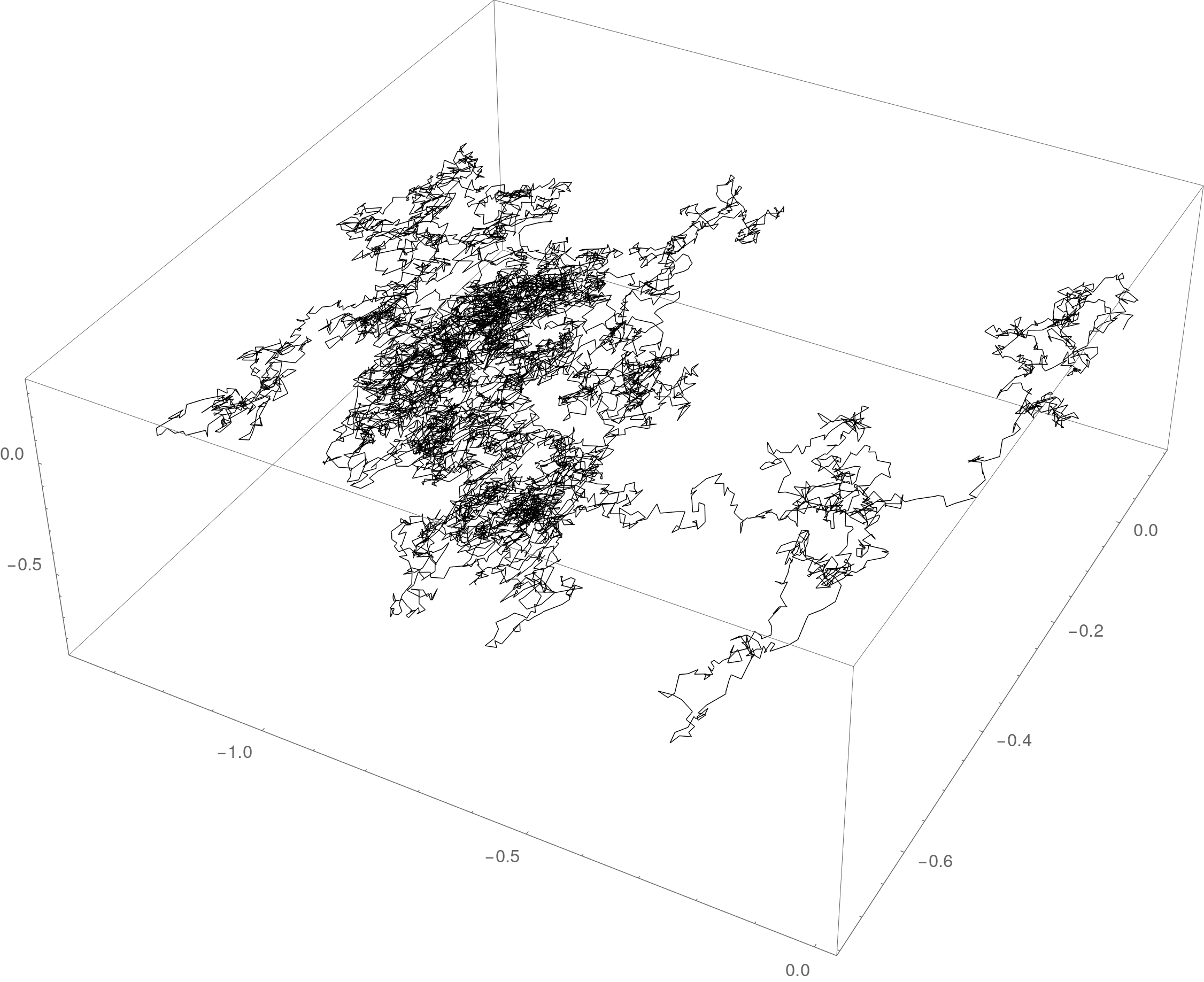}
\hskip-0in\,\\[-.05in]
\caption{$\mathbb{R}^3$ tree process from $3$-ensemble of a planar continuous process starting at $(0,0,0).$}\label{three-state} \vskip-.075in 
\end{figure}
\noindent Such construction can be glued$,$ for all $t\in\mathbb{R}_{\geqslant0},$ in equivalence class (partition) for all $\frac{\eta}{2}\in\mathbb{N}_0,$ by$:$
\\[-.2in] 
\begin{align} 
\mathcal{T}_\eta\,=\,\bigg\{\frac{t}{(\frac{\eta}{2})! ~ 2^{(\eta/2)}}\,: \,~ 0\,\leqslant\,t\,\leqslant f_\eta=\left(\frac{\eta}{2}\right)! ~ 2^{(\eta/2)} \bigg\} 
\hskip.2in \Big| 
\hskip.1in \mathcal{T} = \bigcup^\infty_{\eta=0}\mathcal{T}_\eta
\end{align}
\\[-.175in]
such that $W_i(0)\!=\!0,$ $W_i(1)\!=\!Z^{(i)}_1,$ $\forall\,i\!=\!i,\ldots,\eta;$ $W(t)\!=\!(W_1(t),\ldots,W_\eta(t))^T;$ $Z_t\!=\!(Z^{(1)}_t\!,\ldots,Z^{(\eta)}_t)^T;$ $Z_{t^{k_\eta\,\in\,(1,\,\ldots,\,2^\eta)}_\eta}\!\in\mathcal{N}[0,1];$ and$,$ $\forall\,t^{k_\eta}_\eta\in\mathcal{T}_\eta\backslash\mathcal{T}_{\eta-1},$ $k_\eta\in(1,\,\ldots,\,2^\eta),$
\\[-.2in] 
\begin{align}\hskip-.0in 
W(t^{k_\eta}_\eta) =\, \frac{1}{\eta!}\left(\!W\!\left(\!t^{k_{\eta-1}}_{\eta-1}+\frac{1}{f_{\eta-1}}\right)+W\!\left(\!t^{k_{\eta-1}}_{\eta-1}\right)\!\right) 
\cdot k_\eta\in(1,\,\ldots,\,2^\eta)
+\,\frac{1}{f_\eta}
\,Z_{t^{k_\eta\,\in\,(1,\,\ldots,\,2^\eta)}_\eta}.
\label{eq:interpolation}
\end{align} 
$\,$\\[-.15in]
Clearly$,$ in interpolation$,$ the $\mathbb{R}^{\lvert\mathcal{A}^{\otimes m}\rvert}$ process is given by independent increments$;$ first summand is exponential in addition to the whole pieces being independent affine-random process with respect to normal distribution$;$ the necessary task of showing a {\it threshold of uniform convergence} for this constructed (interpolated) process $W\!:\, [0,1]\longrightarrow\mathbb{R}^{\lvert\mathcal{A}^{\otimes m}\rvert}$ can then follow$.$ \hfill $\square$
\\[-.2in]
\begin{theorem}[meromorphic extension]\label{stationlem} 
On continuous ${\rm\bf P}^{\otimes n}_{\!\!t}$ limit ${\rm\bf P}_{\!(t)},$ Markov process $(X_t)_{t\geqslant0},$
\\[-.15in]
\begin{align} 
\hskip-.1in 
\begin{array}{c} \displaystyle \lim_{n\uparrow\infty}\lvert\mathcal{A}\rvert^{-1}\left\lVert 
{\rm\bf P}^{\otimes n}_{\!\!t} - {\rm\bf P}_{\!(t)} 
\right\rVert_{\mathcal{L}(L^p(\gamma))}\!\!\overset{n\rightarrow\infty}{\longrightarrow}0 
\\[.0in] 
\\[-.15in] 
\hskip.0in \displaystyle \forall ~ p\geqslant 1
\end{array}
\hskip.0in 
\begin{array}{l} 
\\[-.2in]
\left| \hskip-.025in 
\begin{array}{l} 
\\[-.25in]
\displaystyle \lVert X\rVert_{\mathcal{L}(L^p(\gamma))} 
\!:=\big(\mathbb{E}[\lVert X\rVert^p]\big)^{\!1/p}
=\Big(\!\sum_{X<\infty}\hskip-.05in
\big(X^\dagger X\big)^{p/2}\Big)^{\!\!1/p} 
\\[.0in]
\\[-.15in]
\lVert X\rVert_{\mathcal{L}(L^\infty(\gamma))} 
\!:=\sup_{X<\infty} \lVert X\rVert_{\mathcal{L}(L^1(\gamma))}
=\sup_{X<\infty} \lvert X\rvert.
\end{array} \right. 
\end{array} \hskip-.2in
\end{align}
\\[-.15in]
In addition$,$ ${\rm\bf P}_{\!t}$ defines meromorphic extension of locally compact Markov trees$,$ on irreducible Markov chain $(S,{\rm\bf P})$ of periodic transitions$,$ for finite state space $S=(0,1,\ldots,n).$
\end{theorem}
$\,$\\[-.2in] 
{\it Proof.}~~WLOG$,$ set ${\rm\bf P}^{\otimes n}_{\!\!t} \!:= {\rm\bf P}(\frac{1}{\sqrt{n}}\,\widetilde{x}_n\!\leqslant\!x);$ \,$\widetilde{x}_n\!=x_1\!+\cdots+\!x_n;$ ${\rm\bf P}^{(j)}_{\!\!t}\!=\!{\rm\bf P}(x_j\!\leqslant\!x);$ then$,$ for some $f,$ 
\\[-.2in] 
\begin{align}\hskip-.8in 
\frac{1}{\lvert\mathcal{A}\rvert}\left\lVert 
{\rm\bf P}^{\otimes n}_{\!\!t} - {\rm\bf P}_{\!(t)} 
\right\rVert_{\mathcal{L}(L^\infty(\gamma))} := 
\frac{1}{\lvert\mathcal{A}\rvert}\Big( 
\Big\lVert 
\frac{{\rm\bf P}^{(j)}_{\!\!t}}{\sqrt{n}} - {\rm\bf P}_{\!(t)} 
\Big\rVert_{\mathcal{L}(L^\infty(\gamma))} 
\Big)^{\!\!\otimes n} 
\!=\,\sup_n\Big\{\!\frac{1}{n\lvert\mathcal{A}\rvert}\Big|\! 
-\!\!f(n)\,{\rm\bf P}_{\!(t)} 
+\hskip.025in {\rm\bf P}(\widetilde{x}_n\!\leqslant\!x) 
\Big|\Big\} \!\longrightarrow\!0
\end{align}
\\[-.15in] 
a.s$.$ as $n\!\longrightarrow\!\infty;$ similarly$,$ on sets $\lambda,$ $\lambda^{\ast},$ of all eigenvalues of diagonal states ${\rm\bf T},$ resp$.$ ${\rm\bf T}^{\sixptr}.$ 
\\[.05in] 
Let ${\rm\bf P}\!=\!({\rm\bf P}_{\!i,j})_{i,j\geqslant0}$ be transition matrix$,$ defined for ${\rm\bf P}_{\!t},$ where ${\rm\bf P}^n_{\!\!i,j}={\rm\bf P}\{X_{n+m}=j\mid X_m=i\}$ is the probability that the process goes from state $i$ to state $j$ in $n$ transitions$,$ given only temporally homogeneous process having stationary transition probabilities$.$ Moreover$,$ let 
\\[-.25in]
\begin{align}
\theta^{(n)}_{i,i}\,=\,{\rm\bf P}\{X_n=i,\,X_m\neq i,\,\forall\,m\!=\!1,\ldots,n\!-\!1 \,\mid\, X_0=i\}, ~ \forall n\geqslant 1 
\hskip.2in \big| \hskip.1in 
\theta^{(0)}_{i,i} := 0, ~ \forall\,i
\end{align}
\\[-.25in]
be the probability that$,$ starting at $i,$ the first return to state $i$ occurs at $n$th transition$.$ 
\\[.05in]
Clearly$,$
\\[-.4in]
\begin{align}
\theta^{(1)}_{i,i}={\rm\bf P}_{\!i,i}. 
\end{align}
\\[-.275in] 
Considering all possibilities of the process in each mutually exclusive event $E_m,$ $m\!=\!1,\ldots,n,$ given by $X_0\!=\!i,$ $X_n\!=\!i,$ and first return to state $i$ occurring at $m$th transition$:$ Given $\theta^{(m)}_{i,i},$ then for the remaining $n\!-\!m$ transitions$,$ we only have possibilities of the process for which $X_n\!=\!i.$ Hence$,$ by Markov property language$,$ for all $1\leqslant m\leqslant n,$
\\[-.25in]
\begin{align}
{\rm\bf P}\{E_m\} \,=\, {\rm\bf P}\{\text{first return is at $m$th transition $\mid X_0\!=\!i$}\} 
\!\cdot\! {\rm\bf P} \{X_n\!=\!i\mid X_m\!=\!i\}
~ = ~ \theta^{(m)}_{i,i}\,{\rm\bf P}^{n-m}_{\!\!i,i}
\end{align}
\\[-.55in]
\begin{align}
{\rm\bf P}^n_{\!\!i,i}= {\rm\bf P} \{X_n\!=\!i\!\mid\!X_0\!=\!i\} 
= \sum^n_{m=1}\!{\rm\bf P}\{E_m\} 
= \sum^n_{m=1}\!\theta^{(m)}_{i,i}\,{\rm\bf P}^{n-m}_{\!\!i,i}
= \sum^n_{m=0}\!\theta^{(m)}_{i,i}\,{\rm\bf P}^{n-m}_{\!\!i,i} 
\hskip.1in \Big| \, {\rm\bf P}^0_{\!\!i,i}\!=\!1, ~ \theta^{(0)}_{i,i}=0.
\end{align}
\\[-.1in]
Now$,$ define generating function of sequence $\{{\rm\bf P}^n_{\!\!i,j}\},$ resp$.$ $\{\theta^{(n)}_{i,j}\},$ as follows$:$
\\[-.25in]
\begin{align}
\Theta_{i,j}(r) = \sum^\infty_{n=0}r^n\,\theta^{(n)}_{i,j} \hskip.05in, 
\hskip.15in 
P_{i,j}(r) = \sum^\infty_{n=0}r^n\,{\rm\bf P}^n_{\!\!i,j} 
\hskip.2in \Big| \hskip.2in \lvert r\rvert < 1. 
\end{align}
\\[-.3in]
Then by 
\\[-.4in]
\begin{align} 
\sum^\infty_{n=0}r^nP^n = (1-rP)^{-1} = \frac{1}{\det(1-rP)}\Big(\text{Cofactor}\Big(1-rP\Big)\Big)^{\!\rm T}
\end{align}
\\[-.2in]
i.e$.$ on matrix inverse (by transpose of cofactor divided by determinant)$,$ it follows that $P_{i,j}(r)$ can be extended to a meromorphic function on complex plane for all $i,j\in S.$
\\[.05in]
In particular$,$ by
\\[-.2in]
\begin{align} 
\Big(\sum^\infty_{k=0} a_kr^k\Big) 
\Big(\sum^\infty_{\ell=0} b_\ell\,r^\ell\Big)
= \sum^\infty_{n=0}r^n \Big(\sum^\infty_{m=0}a_mb_{n-m}\Big)
=\sum^\infty_{n=0}c_nr^n 
\hskip.2in \Big| \hskip.05in c_n=\sum^n_{m=0}a_mb_{n-m},
\end{align}
\\[-.15in]
we have
\\[-.45in]
\begin{align} 
\Theta_{i,i}(r)\,P_{i,i}(r)\,=\,P_{i,i}(r)-1,
\hskip.1in \text{ i.e$.$ } 
\hskip.1in P_{i,i}(r)\,=\,\frac{1}{1-\Theta_{ii}(r)} 
\hskip.2in \Big| \hskip.2in \lvert r\rvert < 1.
\end{align}
\\[-.15in]
By $\sum^\infty_{n=0}\theta^{(n)}_{i,i}=1,$ then $\lim_{r\longrightarrow 1}(-\Theta_{i,i}(r)) = 1.$ Hence $P_{i,i}(r)$ has a pole at $r=1.$ 
\\[.05in]
As a result$,$ 
\\[-.45in]
\begin{align} 
\Theta_{i,i}(r)\,=\,1-\frac{1}{P_{i,i}(r)}
\end{align}
\\[-.2in]
defines a meromorphic extension of $\Theta_{i,i}(r)$ with removable singularity at $r=1,$ i.e$.$ let $\varrho_{i,i}\in [1,\infty)$ denote radius of convergence of power series for $\Theta_{i,i}(r).$ By $\theta^{(n)}_{i,i}\geqslant 0,$ $\forall n\in\mathbb{N},$ then either $\varrho_{i,i}=\infty$ or $\varrho_{i,i}$ is a singular point for the meromorphic extension of $\Theta_{i,i}(r).$ Thus$,$ $\varrho_{i,i}>1$ for $i\in S.$ \hfill $\square$
\\[.1in]
{\it Remark.}~~Since every local field is either a finite algebraic extension of the $p$-adic number field for prime $p$ or finite algebraic extension of the $p$-series field$;$ moreover$,$ a locally compact$,$ non-discrete$,$ topological field which is not totally disconnected is necessarily either the real or the complex numbers$;$ let the rings of integers for $p$-adic numbers and $p$-series field (i.e$.$ field of formal series with coefficients from finite field of $p$ elements) be represented by a $p$-ary tree (although the $p$-adic field has characteristic $0$ whereas the $p$-series field has characteristic $p$)$,$ then $({\rm\bf P}_{\!t})_{t\geqslant0}$ satisfies family of meromorphic functions of genus $g\!>\!0$ tree having a single pole (an analog of polynomials for higher genera i.e$.$ with weakened planarity condition)$.$ 
\\[.1in]
Thus$,$ we have an object for the flexible classification of meromorphic functions (often called topological classification)$,$ which relates$,$ Fig$.$~\ref{braid-string-link}$,$ tree isometry to enumerative algebraic geometry and singularity theory$,$ on the one hand$,$ and to braid group action on constellations$,$ on the other hand$.$ Classically$,$ in \cite{bb84} Alice$,$ Bob$,$ and eavesdropper Eva known $\{{\rm\bf P}_{\!\!x}\},$ for uniform $\{{\rm\bf P}^{\sixptr}_{\!\!x'}\},$ the total variation distance $\frac{1}{2}\sum_{x\,\in\,\mathcal{X}} \big| {\rm\bf P}_{\!\!x} \!-\! {\rm\bf P}^{\sixptr}_{\!\!x} \big| \!\leqslant\!\varepsilon$ implies the uncorrelated ${\rm\bf P}[x\!\neq\!x']\!\leqslant\!\varepsilon,$ $\forall \,x,x'\!\in\!\mathcal{X}.$ 
$\,$\\[-.1in]
$\,$\\[-.45in]
\begin{definition}\label{def2} 
Let order $\alpha$ divergence of probability density $f$ from density $f_+$ be $h_\alpha(f\lVert f_+)\!:$
\\[-.2in] 
\begin{align}\hskip-.0in 
h_\alpha(f\lVert f_+) = \frac{1}{1-\alpha} 
\ln\!\Bigg(\!\int_{\Omega}
\,\frac{f^\alpha(x)}{f^{\alpha-1}_+(x)}
\,dx\!\Bigg); \hskip.05in  \text{ resp$.$ } \hskip.05in 
\frac{1}{\ln\lvert\mathcal{A}\rvert^{1-\alpha}} \ln\!\Bigg(\sum_{x\,\in\,\mathcal{X}: \,\lvert\mathcal{X}\rvert\,=\,\lvert\mathcal{A}\rvert} 
\,\frac{f^\alpha(x)}{f^{\alpha-1}_+(x)}
\Bigg) 
\end{align}
\\[-.15in]
in the limit $\longrightarrow\alpha,$ for all $\alpha\geqslant0;$ where $\Omega\equiv\mathbb{R}$ resp$.$ $\Omega\equiv\mathbb{R}^n),$ $\forall\,n\in\mathbb{N}.$
\end{definition}
$\,$\\[-.5in]
\begin{lemma}[Ergodic classification]\label{ergodiclem}
Ergodic property of the right continuous process $(X_t)_{t\geqslant0}$ is classified in terms of entropy and divergence$,$ as follows$:$ 
\\[-.225in] 
\begin{align}
&\hskip-.05in h_0(f\lVert f_+)=\, -\ln f_+(\{x\mid f(x)>0\}) 
\hskip.1in \text{ ($-1$ natural log of conditional probability for $f_+$)}
\\[.05in]
&\hskip-.05in h_{1/2}(f\lVert f_+)=\, -2\log\!\int_{\Omega} \sqrt{f(x)f_+(x)}\,dx 
\hskip.1in \text{ ($-2$ natural log of Bhattacharyya coefficient \cite{bh46})}
\\[.025in]
&\hskip-.05in h_1(f\lVert 1)=\, -\!\int_{\Omega} f(x)\ln f(x)\,dx ~ =\, h(f) 
\hskip.1in \text{ (differential entropy)}
\\[.05in]
&\hskip-.05in h_1(f\lVert f_+) =\! \int_{\Omega} f\ln\frac{f}{f_+}\,dx 
\hskip.1in \Big| \hskip.05in 
0\ln\frac{0}{0} := 0 \,\text{ for continuity}
\hskip.1in \text{ (Kullback-Leibler divergence \cite{kl51})}
\\[.05in]
&\hskip-.05in h_2(f\lVert f_+)=\, -\ln \mathbb{E}[f/f_+] 
\hskip.1in \text{ (natural log of expectation of the ratio $f/f_+$)}
\\[.05in]
&\hskip-.05in h_\infty(f\lVert f_+)=\, -\ln\sup_x (f(x)/f_+(x)) 
\hskip.1in \text{ ($-1$ natural log of maximum of the ratio $f/f_+$)}
\\[-.0in]
&\hskip-.05in h_\alpha(f\lVert 1)=\, \frac{1}{1-\alpha} \ln\bigg(\!\int_{\Omega}\hskip-.05in f^\alpha(x)
\,dx\!\!\bigg), \hskip.0in  \text{ resp$.$ } \hskip.05in 
\frac{1}{\ln\lvert\mathcal{A}\rvert^{1-\alpha}} \ln\!\Bigg(\!\sum_{x\,\in\,\mathcal{X}: \,\lvert\mathcal{X}\rvert\,=\,\lvert\mathcal{A}\rvert} 
\hskip-.3in f^\alpha(x)\!\Bigg) 
\hskip.0in \text{ (modified R\'enyi \cite{re61})}.
\end{align}
\end{lemma}
$\,$\\[-.1in]
{\it Proof.}~~The proof follows by taking respective limit on the definition$.$ \hfill $\square$
\\[.1in]
{\it Remark.}~~The discrete equivalence with $\sum$ replaces $\int$ is R\'enyi (divergence) entropy classification$.$ 
\\[-.2in] 
\begin{theorem}[asymptotic equipartition property] 
Let $X_1,\ldots,X_n$ be a sequence of random variables drawn i.i.d$.$ according to a probability measure $\mu$ density $f(x).$ Then 
\\[-.25in] 
\begin{align}
-\,\frac{1}{n}\ln f^{\otimes n} (X_1,\ldots,X_n)
\hskip.05in\longrightarrow\hskip.05in 
\mathbb{E}[-\ln f(X)]\,=h_1(f\lVert 1)=\,h(f)
\end{align}
\\[-.225in]
almost surely$,$ where $f^{\otimes n}$ denotes the $n$th convolution$.$ 
\end{theorem}
$\,$\\[-.25in]
{\it Proof.}~~This follows directly from the strong law of large numbers$,$ for all $f$ with respect to $\mu.$ \hfill $\square$ 
\\[.1in] 
{\it Remark.}~~Ergodicity$,$ Fig$.$~\ref{2dentropy}$,$ defined in terms of $1$-divergence (differential entropy)$,$ exists in limit$:$
\\[-.25in] 
\begin{align}
h(f)=h_1(f\lVert 1)= \lim_{n\longrightarrow\infty} 
\frac{1}{n} h(f^{\otimes n}) 
\hskip.1in \Big| \hskip.1in 
h(f^{\otimes n}) = -\ln f^{\otimes n} (X_1,\ldots,X_n).
\end{align}
\vskip-.2in $\,$
$,$\\[-.45in]
\begin{definition}\label{def3} 
Let $\mu$ be continuous probability measure on $\mathcal{X}.$ We say $\rho$ is $\mu$-entire$,$ if 
\\[-.2in] 
\begin{align}
\rho(z) = e^{P_n(z)} = \int_{\Omega}\widetilde{\rho}(zx)\mu(dx) 
\hskip.2in \Big| \hskip.05in \Omega\equiv\mathbb{R} \text{ $($resp$.$ } \Omega\equiv\mathbb{R}^n)
\end{align}
\\[-.15in]
holds for all $z\in\mathbb{C}$ $($resp$.$ $z\in\mathbb{C}^n),$ and degree $n\!\geqslant\!1$ polynomial $P_n$ with complex coefficients$.$
\end{definition}
$\,$\\[-.2in]
It is not a difficult task to show that all bounded $\mu$-entire functions on $\mathcal{X}$ are constant$.$ Thus$,$ for $h(f)<\infty,$ the vanishing of entropy $h$ is equivalent to trivialness of Poisson boundary $\Pi(\mathcal{X},f)$ \cite{kv83}$.$
\\[.1in]
{\it Remark (Hadamard's factorization theorem).}~~Let $\rho(z)$ be  entire function$.$ Assume for simplicity that $\rho(0)\neq0.$ Then$,$ the {\it genus} of an entire function is the smallest integer $h$ such that $\rho(z)$ can be represented in the form
\\[-.3in]
\begin{align}
\rho(z) = e^{\theta(z)}\prod_n\Big(1-\frac{z}{a_n}\Big) 
\exp\Big(z/a_n + \cdots + \frac{1}{2}(z/a_m)^m\Big)
\end{align}
\\[-.15in]
where $\theta(z)$ is degree $\leqslant\! m$ polynomial$.$ And$,$ if there is no such representation$,$ the genus is infinite$.$
\\[.05in]
Denoting by $M(r)$ the maximum of $\lvert\rho(z)\rvert$ on $\lvert z\rvert=r.$ The order of the function $\rho(z)$ is 
\\[-.2in]
\begin{align} 
\lambda = \limsup_{r\longrightarrow\infty} \frac{\ln\ln M(r)}{\ln r}.
\end{align}
\\[-.15in]
In addition$,$ according to this definition$,$ $\lambda$ is the smallest number such that 
\\[-.2in]
\begin{align}
M(r) \leqslant e^{r^{\lambda\,+\,\varepsilon}}
\end{align}
\\[-.25in]
for any given $\varepsilon>0,$ as soon as $r$ is sufficiently large$.$ Moreover$,$ it is known (theorem) that the genus and the order of an entire function satisfy
\\[-.25in]
\begin{align}
h \leqslant \lambda \leqslant h+1.
\end{align}
$\,$\\[-.5in]
\subsection*{Classical (diagonal) operator metrizing for compact uniform Markov trees} 

\begin{lemma}[classical bipartite]\label{class-bipartite-lemma} 
Suppose ${\rm\bf P}_{\!X}(x)\!\mid\!x\!\in\!\mathcal{X},$ compact separable Hilbert $L^p(\gamma)$ bipartite space $\mathcal{X}_K\!\oplus\!\mathcal{X}_G,$ for $\lvert\xi\rvert\!=\!n\!+\!1$ mixture$,$ with unique mixed state$,$ on stopping time $t$ process $(X_t)_{t\geqslant0}$ where $K\!=\!G(X),$ for uniform-density family $G$ of $\lvert\mathcal{A}\rvert$-universal hash functions i.e$.$ $\kappa\!\in\!K,\,\forall\,g,$ by
\\[-.25in]
\begin{align} 
G \ni g\!: \,\mathcal{A}^{\otimes \displaystyle m}\cong\mathcal{X} \longrightarrow \mathcal{A}^{\otimes \displaystyle k} \hskip.075in \mid 
\hskip.075in \mathcal{A}^{\otimes \displaystyle \zeta}\,\cong\,G. 
\end{align}
\\[-.275in]
For $\alpha\!=\!1,$ minimum entropy $h_{\min};$ joint conditional density ${\rm\bf T}_{\!KG},$ resp$.$ uniform ${\rm\bf T}^{\sixptr}_{\!KG},$ where $p\!=\!1,$ 
\\[-.225in]
\begin{align} 
\hskip-.0in\lvert\mathcal{A}\rvert^{-1} 
\left\lVert{\rm\bf T}_{\!KG}-{\rm\bf T}^{\sixptr}_{\!KG}\right\rVert_{\mathcal{L}(L^p(\gamma))} 
~ \leqslant \,~  \lvert\mathcal{A}\rvert^{^{\displaystyle\!\Big(\!\!-\frac{h_+\!-\!k}{2}\Big)}}\,, 
\hskip.2in \forall ~ h_{\min}\geqslant h_+.
\end{align}\\[-.5in]
\end{lemma}
$\,$\\[-.175in] 
{\it Proof.}~~Assuming discrete WLOG$,$ for $\zeta\!\neq\!k$ on $G$ uniform$:$ ${\rm\bf P}_{\!G}(g)=\lvert\mathcal{A}\rvert^{\displaystyle-\zeta}.$ For $\lvert K\rvert\!=\!\lvert\mathcal{A}\rvert^{\displaystyle k}$ by $G$ universal$,$ then the uniform ${\rm\bf P}^{\sixptr}_{\!KG}(\kappa,g)$ (by independence) and conditional ${\rm\bf P}_{\!KG}(\kappa,g)\!:$
\\[-.225in]
\begin{align}
{\rm\bf P}^{\sixptr}_{\!KG}(\kappa,g)\,=\,
\lvert\mathcal{A}\rvert^{\displaystyle-k\!-\!\zeta}\,, 
\hskip.25in 
{\rm\bf P}_{\!KG}(\kappa,g)\,=\,{\rm\bf P}_{\!G}(g)\,{\rm\bf P}_{\!\{K|G\}}(\kappa|g)\,=\,
\lvert\mathcal{A}\rvert^{\displaystyle-\zeta}\hskip-.15in \sum_{x\,:\,\,\kappa\,=\,g(x)}\hskip-.15in 
{\rm\bf P}_{\!\!x}
\end{align}
\\[-.15in]
where the RHS sum is ${\rm\bf P}(K\!\!=\!\kappa),$ admitting collision principle$:$ for $K$ small enough$,$ $K$ and $G$ joint state density gets close to uniform (independence$,$ no-correlation)$,$ in order for no information to be revealed by some eavesdropped hash function $g$ in an open receiver-sender channel$.$ 
\\[.075in]
Following existence$,$ and for all central limit of subsequences$,$ by the trace norm
\\[-.25in]
\begin{align} 
\left\lVert{\rm\bf T}_{\!KG} - {\rm\bf T}^{\sixptr}_{\!KG} \right\rVert_{\mathcal{L}(L^1(\gamma))} 
\,=~ 
\sum_{\kappa,\,g}\Big| 
{\rm\bf P}_{\!\!\kappa,g} - 
\lvert\mathcal{A}\rvert^{\displaystyle-k\!-\!\zeta}\Big|
\nonumber 
\end{align}
\\[-.4in] 
then 
\\[-.275in]
\begin{align}\hskip-.1in 
\lvert\mathcal{A}\rvert^{-1}\left\lVert{\rm\bf T}_{\!KG} - {\rm\bf T}^{\sixptr}_{\!KG} \right\rVert_{\mathcal{L}(L^1(\gamma))} 
~=~ 
\lvert\mathcal{A}\rvert^{-1}\sum_{\kappa,\,g} 
\bigg|\lvert\mathcal{A}\rvert^{\displaystyle-\zeta}\hskip-.175in 
\sum_{x\,:\,\,\kappa\,=\,g(x)} 
\hskip-.175in {\rm\bf P}_{\!\!x} - 
\lvert\mathcal{A}\rvert^{\displaystyle-k\!-\!\zeta}\bigg| 
~=~ 
\lvert\mathcal{A}\rvert^{\displaystyle-\zeta-1}\sum_{\kappa,\,g} \bigg|\sum_{x\,:\,\,\kappa\,=\,g(x)} 
\hskip-.175in {\rm\bf P}_{\!\!x} - 
\lvert\mathcal{A}\rvert^{\displaystyle-k}\bigg| 
~ \doteq  
\nonumber
\end{align}
\\[-.15in]
i.e$.,$ by Cauchy-Bunyakovsky-Schwarz inequality 
\\[-.225in]
\begin{align} 
\bigg|\sum_{x\,\in\,\mathcal{X}}\!a_x\,b_x\bigg| 
\,\leqslant\,\bigg(\sum_{x\,\in\,\mathcal{X}} \!a^2_x\!\bigg)^{\displaystyle \!1/2} \hskip.025in 
\bigg(\sum_{x\,\in\,\mathcal{X}}
\!b^2_x\!\bigg)^{\displaystyle \!1/2}; 
\hskip.2in 
\bigg|\sum_{x\,\in\,\mathcal{X}}\!a_x\bigg| 
\,\leqslant\,\bigg(\sum_{x\,\in\,\mathcal{X}} \!a^2_x\!\bigg)^{\displaystyle \!1/2} \hskip.025in 
\bigg(\sum_{x\,\in\,\mathcal{X}}1\!\bigg)^{\displaystyle \!1/2} \,=\,
\sqrt{|X|}~\bigg(\sum_{x\,\in\,\mathcal{X}} 
\!a_x^2\!\bigg)^{\displaystyle \!1/2};
\nonumber
\end{align}
\\[-.4in]
\begin{align} 
\hskip-.7in \doteq ~ 
\lvert\mathcal{A}\rvert^{\displaystyle-\zeta-1}\cdot 
\lvert\mathcal{A}\rvert^{\displaystyle ((k\!+\!\zeta)/2)}
\Bigg(\!\sum_{\kappa,\,g}\bigg| 
\sum_{x\,:\,\,\kappa\,=\,g(x)} 
\hskip-.175in {\rm\bf P}_{\!\!x} - 
\lvert\mathcal{A}\rvert^{\displaystyle-k}\bigg|^{\displaystyle 2} 
\,\Bigg)^{\displaystyle \!1/2} 
= ~ 
\lvert\mathcal{A}\rvert^{\displaystyle ((k/2)\!-\!1)} 
\Bigg(\!\lvert\mathcal{A}\rvert^{\displaystyle-\zeta}\,\sum_{\kappa,\,g}\bigg| 
\sum_{x\,:\,\,\kappa\,=\,g(x)} 
\hskip-.175in {\rm\bf P}_{\!\!x} - 
\lvert\mathcal{A}\rvert^{\displaystyle-k}\bigg|^{\displaystyle 2} 
\,\Bigg)^{\displaystyle \!1/2} 
\nonumber
\end{align}
\\[-.45in]
\begin{align} 
\hskip-.8in = ~ 
\lvert\mathcal{A}\rvert^{\displaystyle ((k/2)\!-\!1)} 
\bigg(\!\lvert\mathcal{A}\rvert^{\displaystyle \zeta}\sum_{\kappa,\,g}\Big|
{\rm\bf P}_{\!\!\kappa,g} - \lvert\mathcal{A}\rvert^{\displaystyle -k\!-\!\zeta} 
\Big|^{\displaystyle 2}\bigg)^{\displaystyle \!1/2} 
= ~ 
\lvert\mathcal{A}\rvert^{\displaystyle ((k/2)\!-\!1)} 
\bigg(\!\lvert\mathcal{A}\rvert^{\displaystyle \zeta}\sum_{\kappa,\,g} 
\!\bigg(\!{\rm\bf P}^2_{\!\!\kappa,g}-2\,{\rm\bf P}_{\!\!\kappa,g}\cdot 
\lvert\mathcal{A}\rvert^{\displaystyle -k\!-\!\zeta} + \lvert\mathcal{A}\rvert^{\displaystyle-2(k\!+\!\zeta)} 
\!\bigg)\!\bigg)^{\displaystyle \!1/2} 
\nonumber 
\end{align}
\\[-.45in]
\begin{align} 
\hskip-.85in =~ \lvert\mathcal{A}\rvert^{\displaystyle ((k/2)\!-\!1)} 
\bigg(\!\lvert\mathcal{A}\rvert^{\displaystyle \zeta} 
\bigg(\!\sum_{\kappa,\,g}{\rm\bf P}^2_{\!\!\kappa,g}\,-\,2\cdot \lvert\mathcal{A}\rvert^{\displaystyle -k\!-\!\zeta}\,+\, 
\lvert\mathcal{A}\rvert^{\displaystyle -2(k\!+\!\zeta)\!+\!k\!+\!\zeta}\!\bigg) 
\!\bigg)^{\displaystyle \!1/2} 
\!=\,~
\lvert\mathcal{A}\rvert^{\displaystyle ((k/2)\!-\!1)} 
\bigg(\!\lvert\mathcal{A}\rvert^{\displaystyle \zeta}\sum_{\kappa,\,g} 
{\rm\bf P}^2_{\!\!\kappa,g} - 
\lvert\mathcal{A}\rvert^{\displaystyle -k} 
\!\bigg)^{\displaystyle \!1/2}\!. \hskip.1in (\ast)
\nonumber
\end{align}
\\[-.15in]
In addition$,$
\\[-.35in]
\begin{align} 
\hskip-.9in\sum_{\kappa,\,g}{\rm\bf P}^2_{\!\!\kappa,g} 
\,&= \sum_{\kappa,\,g}
\,\lvert\mathcal{A}\rvert^{\displaystyle-2\zeta}
\hskip-.275in \sum_{\substack{x,\,x':\\[.025in] 
\kappa\,=\,g(x)\,=\,g(x')}} 
\hskip-.275in {\rm\bf P}_{\!\!x}\,{\rm\bf P}_{\!\!x'}
\,=\,\lvert\mathcal{A}\rvert^{\displaystyle-2\zeta} 
\!\left(\!\sum_{\kappa,\,g}
\!\left(\hskip-.075in 
\begin{array}{l} 
\displaystyle \sum_{\substack{x\,\neq\,x':\\[.025in] 
\kappa\,=\,g(x)\,=\,g(x')}} 
\hskip-.275in {\rm\bf P}_{\!\!x}\,{\rm\bf P}_{\!\!x'} 
\,~ + \hskip-.2in 
\sum_{\substack{x\,=\,x':\\[.05in] 
\kappa\,=\,g(x)\,=\,g(x')}}\hskip-.275in 
{\rm\bf P}^2_{\!\!x} 
\end{array}\right)\!\right) 
\!=\,\lvert\mathcal{A}\rvert^{\displaystyle-k\!-\!\zeta} 
\!\sum_{x\,\neq\,x'} \!{\rm\bf P}_{\!\!x}\,{\rm\bf P}_{\!\!x'} 
\,+\,\lvert\mathcal{A}\rvert^{\displaystyle-\zeta}\sum_x{\rm\bf P}^2_{\!\!x} \hskip.1in 
\text{\circled{$\leqslant$}} 
\nonumber
\end{align}
\\[-.1in]
for $\kappa$-partition$,$ $\sum_{\kappa}\sum_{x:\,\kappa\,=\,g(x)} \!=\! \sum_x,$ on ${\rm\bf P}[G(x)\!=\!G(x')]\leqslant\lvert\mathcal{A}\rvert^{\displaystyle -k},$ \,$\lvert\{g\!:g(x)\!=\!g(x')\}\rvert\leqslant\lvert\mathcal{A}\rvert^{\displaystyle \zeta\!-\!k}.$ 
\\[.05in]
Moreover$,$ the sum $\sum_{\kappa,\,g}{\rm\bf P}^2_{\!\!\kappa,g},$ called collision probability$,$ is bounded by $\lvert\mathcal{A}\rvert^{\displaystyle -k\!-\!\zeta}$ as follows$:$ 
\\[-.15in]
\begin{align} 
\sum_{\kappa,\,g}\!\sum_{\substack{x\,\neq\,x':\\[.025in] 
\kappa\,=\,g(x)\,=\,g(x')}} 
\hskip-.275in {\rm\bf P}_{\!\!x}\,{\rm\bf P}_{\!\!x'} 
\hskip.1in = \hskip.1in 
\sum_g\!\!\sum_{\substack{x\,\neq\,x':\\[.025in] 
g(x)\,=\,g(x')}}
\hskip-.175in {\rm\bf P}_{\!\!x}\,{\rm\bf P}_{\!\!x'} 
\hskip.1in = \hskip.1in 
\sum_{x\,\neq\,x'} \,~ \sum_{g(x)\,=\,g(x')} 
\hskip-.175in {\rm\bf P}_{\!\!x}\,{\rm\bf P}_{\!\!x'} 
\hskip.1in = \hskip.1in 
\sum_{x\,\neq\,x'}\!{\rm\bf P}_{\!\!x}\,{\rm\bf P}_{\!\!x'} 
\hskip-.15in \sum_{g(x)\,=\,g(x')}\hskip-.1in 1 
\nonumber 
\end{align}
\\[-.4in]
\begin{align} 
\leqslant \hskip.05in \sum_{x\,\neq\,x'}
{\rm\bf P}_{\!\!x}\,{\rm\bf P}_{\!\!x'}\cdot 
\lvert\mathcal{A}\rvert^{\displaystyle \zeta\!-\!k} 
\hskip.1in \leqslant \hskip.1in 
\sum_{x,\,x'}\,{\rm\bf P}_{\!\!x}\,{\rm\bf P}_{\!\!x'}\cdot 
\lvert\mathcal{A}\rvert^{\displaystyle \zeta\!-\!k} 
\hskip.1in = \hskip.1in 
\lvert\mathcal{A}\rvert^{\displaystyle \zeta\!-\!k} 
\nonumber 
\end{align}
\\[-.125in]
where the last inequality stems from all i.i.d $\sum^{m-1}_{i=0}{\rm\bf P}_{\!\xi_i}=1 \mid {\rm\bf P}_{\!\xi_0\cdots\,\xi_{m-1}}\!\!=\prod^{m-1}_{i=0}{\rm\bf P}_{\!\xi_i},$  
\\[-.225in]
\begin{align} \hskip-.0in 
\sum_{\xi_0,\ldots,\,\xi_{m-1}} 
~ \prod^{m-1}_{i=0}{\rm\bf P}_{\!\xi_i} ~ = ~  
\sum_{\xi_0=\,\cdots\,=\,\xi_{m-1}} 
~ \prod^{m-1}_{i=0}{\rm\bf P}_{\!\xi_i} ~ + ~  
\sum^{m-1}_{j=1} \hskip.1in (j+1)! \hskip-.5in 
\sum_{\xi_0<\,\cdots\,<\,\xi_j;
~ \xi_{j+1}=\,\cdots\,=\,\xi_{m-1}} 
~ \prod^{m-1}_{i=0}{\rm\bf P}_{\!\xi_i}
\nonumber
\end{align}
\\[-.125in]
which equals $1$ if \,$\sum^{m-1}_{i=0}{\rm\bf P}_{\!\xi_i}=1;$ moreover$,$ for all random graph ${\rm sgn}(\pi_\xi\!\mid\!\xi\!=\!(\xi_0,\ldots,\xi_{m-1}))\!:$
\\[-.225in]
\begin{align} \hskip-.65in 
\sum_{\xi_0=\,\cdots\,=\,\xi_{m-1}} \hskip-.05in
\bigg(\prod^{m-1}_{i=0}{\rm\bf P}_{\!\xi_i}\!\bigg)
\,{\rm sgn}(\pi_\xi)\,a_{\xi_0}\!\otimes\cdots\otimes a_{\xi_{m-1}} ~ + ~ \sum^{m-1}_{j=1} \hskip.1in (j+1)! \hskip-.5in 
\sum_{\xi_0<\,\cdots\,<\,\xi_j;
~ \xi_{j+1}=\,\cdots\,=\,\xi_{m-1}} \hskip-.05in
\bigg(\prod^{m-1}_{i=0}{\rm\bf P}_{\!\xi_i}\!\bigg)
\,{\rm sgn}(\pi_\xi)\,a_{\xi_0}\!\otimes\cdots\otimes a_{\xi_{m-1}}.
\nonumber
\end{align}
\\[-.125in]
Then$,$ following the second-order modified R\'{e}nyi entropy$,$ 
\\[-.225in]
\begin{align} 
\sum_x{\rm\bf P}^2_{\!\!x} 
\,=\lvert\mathcal{A}\rvert^{\displaystyle-h_2(x)} 
\leqslant \lvert\mathcal{A}\rvert^{\displaystyle-h_+} 
\hskip.1in \bigg| \hskip.05in 
h_2(x)\,=\, 
\frac{1}{-\ln\lvert\mathcal{A}\rvert} \ln\!\Bigg(\sum_{x\,\in\,\mathcal{X}: \,\lvert\mathcal{X}\rvert\,=\,\lvert\mathcal{A}\rvert} 
\hskip-.25in {\rm\bf P}^2_{\!\!x}\hskip.075in\Bigg) 
\geqslant\,h(x)\,\geqslant\,h_{\min}\,\geqslant\,h_+. 
\nonumber
\end{align}
\\[-.15in]
That is$,$ the assertion follows$:$
\\[-.4in]
\begin{align}
\text{\circled{$\leqslant$}} \hskip.1in 
\lvert\mathcal{A}\rvert^{\displaystyle-k\!-\!\zeta} \,+\, 
\lvert\mathcal{A}\rvert^{\displaystyle-h_+\!-\!\zeta}\,. 
\nonumber
\end{align}
\\[-.25in]
Hence$,$ $(\ast)$ implies$:$ 
\\[-.4in]
\begin{align} 
\lvert\mathcal{A}\rvert^{-1}\left\lVert{\rm\bf T}_{\!KG} - {\rm\bf T}^{\sixptr}_{\!KG} 
\right\rVert \,\,\leqslant\,\, \lvert\mathcal{A}\rvert^{^{\hskip-.025in \Big(\displaystyle\frac{k}{2}\!-\!1\!\Big)}} \!\Big(\lvert\mathcal{A}\rvert^{\displaystyle-k} 
+\,\lvert\mathcal{A}\rvert^{\displaystyle-h_+} 
-\,\lvert\mathcal{A}\rvert^{\displaystyle-k}\Big)^{\!1/2} 
=\,\,\lvert\mathcal{A}\rvert^{^{\hskip-.025in \Big(\displaystyle\frac{k\!-\!h_+}{2}-1\!\Big)}} 
<\,\, \lvert\mathcal{A}\rvert^{^{\hskip-.025in \Big(\displaystyle\frac{k\!-\!h_+}{2}\Big)}}.
\nonumber
\end{align}
\\[-.175in] 
Hence$,$ the conclusion$.$ \hfill $\square$
\\[.1in]
{\it Remark.}~~The Von Neumann uncertainty dualities to estimate entropy $h_{\min}$ is given as follows$:$
\\[-.25in]
\begin{align} 
&\hskip-.0in h_{\min}(A|B) + h_{\max}(A|C) = 0; 
\hskip.05in h_\alpha(A|B) + h_\beta(A|C) = 0, 
\hskip.05in (\alpha\beta)^{-1}(\alpha\!+\!\beta)=2; 
\hskip.05in h_{\frac{1}{2}}\!\leqslant\!h_0; 
\hskip.05in h_{\frac{2}{3}}\!\leqslant\!h_0 
\\[.025in] 
&\hskip-.0in h_{\min}(X|B) + h_{\max}(Z|C) 
\geqslant \log(c^{-1}); 
\hskip.15in h_2(X|B) + h_{\frac{2}{3}}(Z|C) 
\geqslant \log(c^{-1}); 
\\[.0in]
&\hskip-.0in h_{\max}(Z|C) \leqslant h_0(Z|C); 
\hskip.15in h_{\frac{2}{3}}(Z|C) \leqslant h_0(Z|C).
\end{align}
\\[-.475in]
\begin{theorem}[classical bipartite]\label{class-bipartite-theorem} 
Following the prior lemma \ref{class-bipartite-lemma}$,$ for $p\!=\!1,$ 
\\[-.2in]
\begin{align} 
\hskip-.0in\lim_{n_t\uparrow\infty}(n_t\lvert\mathcal{A}\rvert)^{-1}\big(\!\left\lVert{\rm\bf T}_{\!KG}-{\rm\bf T}^{\sixptr}_{\!KG}\right\rVert_{\mathcal{L}(L^p(\gamma))} 
\!\big)^{\!\otimes n_t}  
~ \leqslant \,~  \lvert\mathcal{A}\rvert^{^{\displaystyle\!\Big(\!\!-\frac{h_+\!-\!k}{2}\Big)}}\,, 
\hskip.2in \forall ~ h_{\min}\geqslant h_+.
\end{align}
\end{theorem}
$\,$\\[-.15in] 
{\it Proof.}~~The proof follows from the prior lemma in strong law of large numbers$.$ \hfill $\square$
\\[-.2in]
\begin{corollary}[diagonal bipartite]\label{class-bipartite-corollary} 
The maximum extractable key length
\\[-.225in]
\begin{align} 
\max(k) \,=\, 
\Big\lfloor 
h_+-\,\frac{2}{\ln\lvert\mathcal{A}\rvert}\ln\!\Big(\lvert\mathcal{A}\rvert^{-1}\lim_{n_t\uparrow\infty}(n_t)^{-1}\big(\!\left\lVert{\rm\bf T}_{\!KG}-{\rm\bf T}^{\sixptr}_{\!KG}\right\rVert_{\mathcal{L}(L^p(\gamma))} 
\!\big)^{\!\otimes n_t}\Big)\Big\rfloor\,, 
\hskip.2in \forall ~ h_{\min}\geqslant h_+.
\nonumber
\end{align}
\\[-.15in]
for all tree-valued process $h_{\min}$ low-enough metric-strong protocol of indistinguishable $\kappa$ in reality$.$
$\,$\\[.1in]
{\rm {\it Proof.}~~The proof follows directly from the prior theorem \ref{class-bipartite-theorem}$.$} 
\end{corollary} 
$\,$\\[-.5in]
\subsection*{Quantum operator metrizing for compact uniform Markov trees} $\,$\vskip-.2in

\begin{figure}[H]
\centering ~ \hskip-6.025in 
\includegraphics[width=.3\textwidth,angle=0,origin=c]
{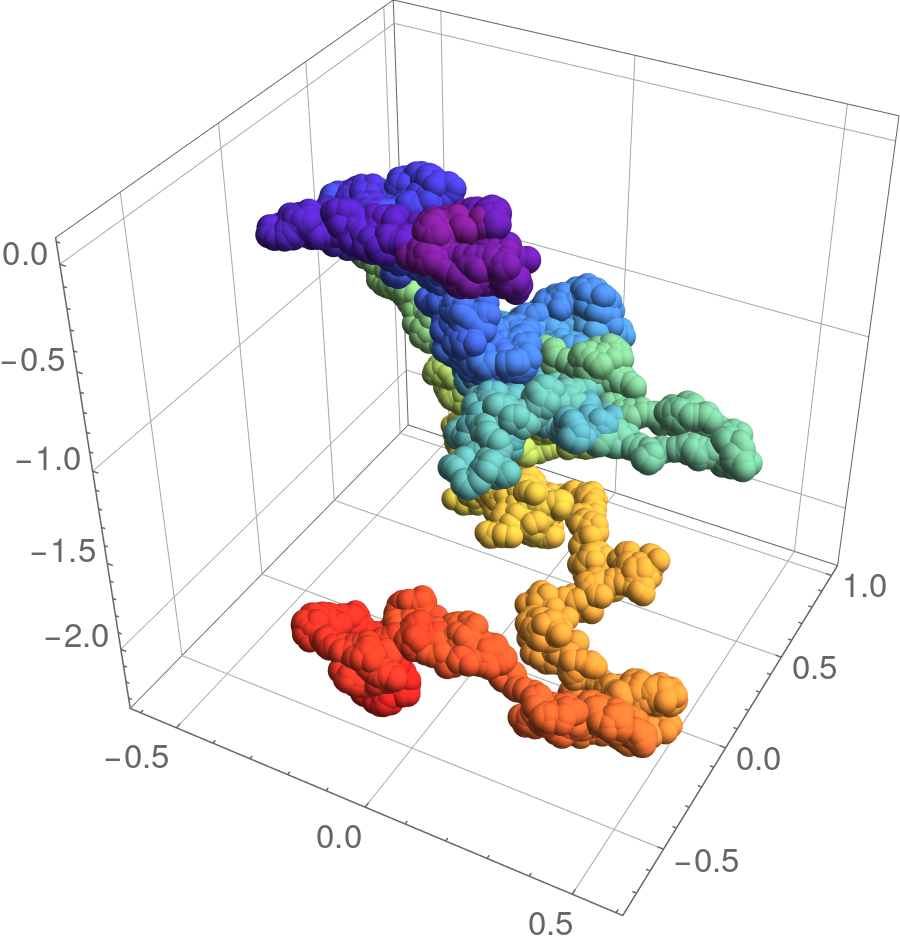} ~ ~ \includegraphics[width=.3\textwidth,angle=0,origin=c]
{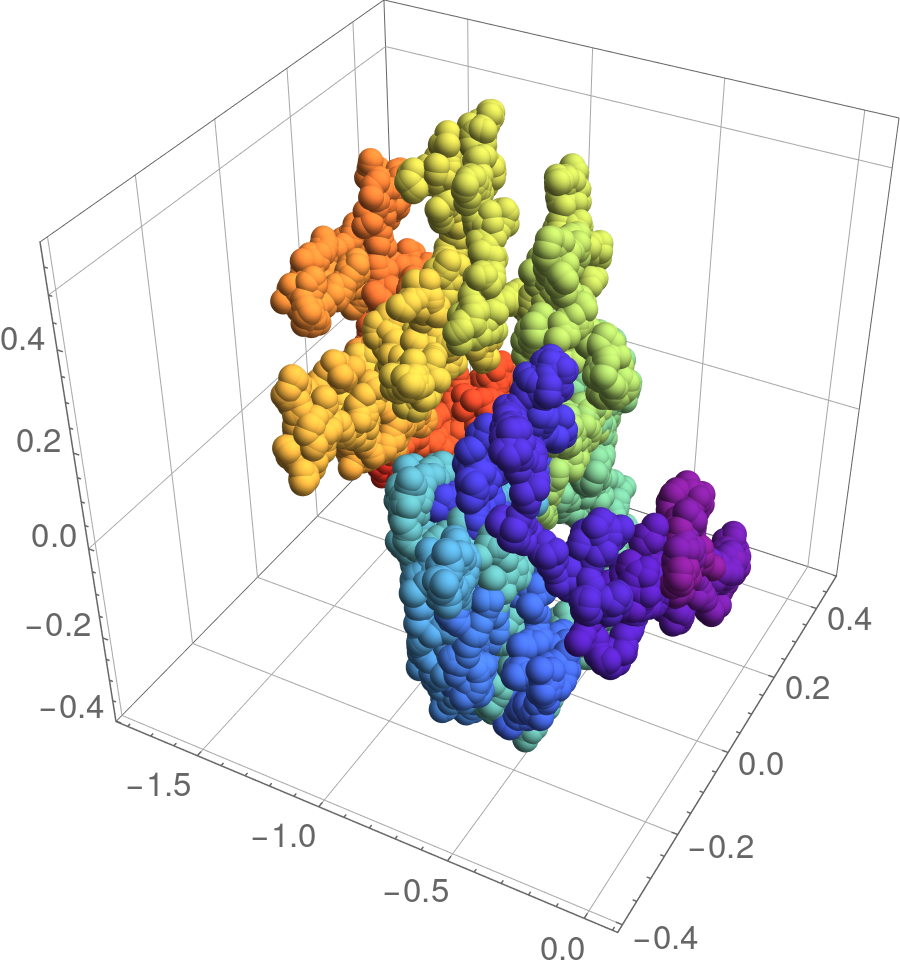} ~ ~ \includegraphics[width=.3\textwidth,angle=0,origin=c]
{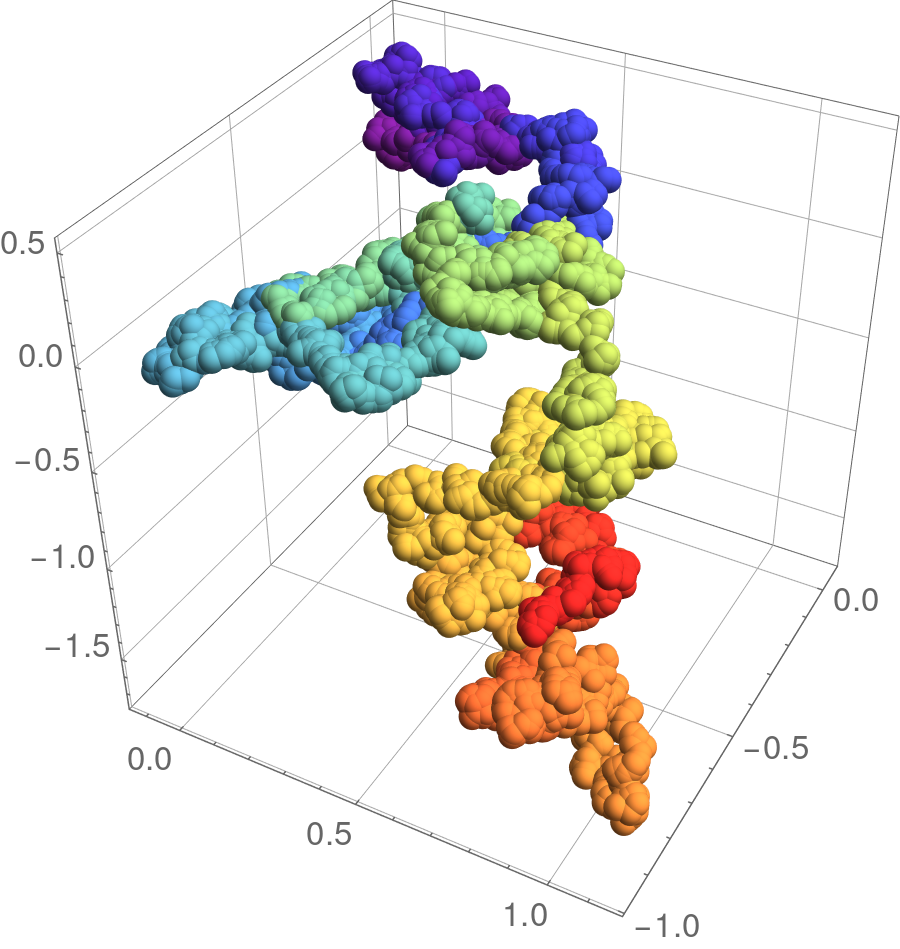}
\hskip-6in\,\\[-.075in]
\caption{$\mathbb{R}^3$ tree-state processes of pairwise keys$;$ partition $\cong$ ring of blocks of Hopf fibration on torus $\mathbb{S}^2.$}\label{hopf-fibration-torus-treeprocess}
\end{figure}

\begin{figure}[H]
\centering ~ \hskip-0.in 
\includegraphics[width=.45\textwidth,angle=0,origin=c]
{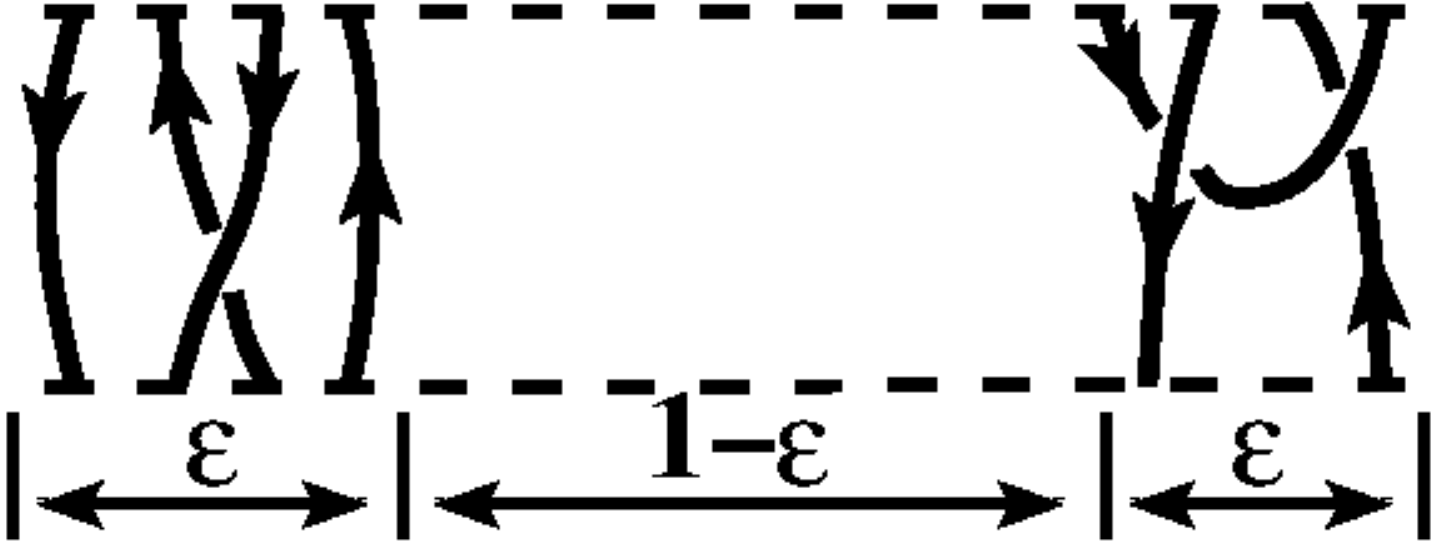}
\hskip-0in\,\\[-.05in]
\caption{$\varepsilon$-rescaling of $\mathbb{R}^3\!:$ $(z,t)\!\longmapsto\!(\varepsilon z,t)$ on $\varepsilon$-parameterized tree tangle $T_1\!\otimes_\varepsilon\!T_2;$ inter-tangle distance $1\!-\!\varepsilon.$}\label{tangle-rescale} \vskip-.075in 
\end{figure}

\noindent The supersymmetric tree process can be encoded by $(\mathcal{P},\gamma)$ bilinear form $\Phi$ in pairwise keys$.$ Fig$.$~\ref{hopf-fibration-torus-treeprocess} does it exactly$:$ $\Phi$ is measure on pairwise ring-blocks of Hopf fibration on torus $\mathbb{S}^2\!:$ in stereographic projection of $\mathbb{S}^3$ to $\mathbb{R}^3$ and compression of $\mathbb{R}^3$ to torus $\mathbb{S}^2$ boundary$,$ following embedding of knot fiber (space $\mathbb{S}^1,$ circle) bundle in total space ($\mathbb{S}^3,$ $3$-sphere)$;$ for properties of trace distance$:$
\\[-.225in]
\begin{align} 
\lvert\mathcal{A}\rvert^{-1}\left\lVert 
\Phi({\rm\bf T}) - \Phi({\rm\bf T}^{\sixptr}) 
\right\rVert_{\mathcal{L}(L^1(\gamma))} 
\leqslant 
\lvert\mathcal{A}\rvert^{-1}\left\lVert {\rm\bf T} - {\rm\bf T}^{\sixptr}\right\rVert_{\mathcal{L}(L^1(\gamma))} 
\leqslant 
\varepsilon 
\end{align} 
\\[-.475in] 
\begin{align}
\lvert\mathcal{A}\rvert^{-1}\,tr\big(\left\lvert 
\Phi({\rm\bf T}) - \Phi({\rm\bf T}^{\sixptr}) 
\right\rvert\big) = 
\lvert\mathcal{A}\rvert^{-1}\,tr\big(\left\lvert 
\Phi({\rm\bf T}-{\rm\bf T}^{\sixptr}) 
\right\rvert\big); 
\hskip.1in 
tr(\Phi({\rm\bf T}))=tr({\rm\bf T}); 
\hskip.1in 
tr(\left\lvert\Phi({\rm\bf T})\right\rvert) 
\neq tr(\left\lvert{\rm\bf T}\right\rvert)
\nonumber 
\end{align} 
\\[-.225in] 
i.e$.$ on \cite{cgob19,bgiko19} formalism of real-tree random metric \cite{defkz00} through $\Phi$-channel quantum-ensemble 
\\[-.25in]
\begin{align}
\big\{{\rm\bf P}_{\!\!x},\widetilde{\rm\bf T}_{\!x}\big\}_{x\,\in\,\mathcal{X}} ~ ; \hskip.1in
\forall\,x\equiv x_{\xi_0}\!\cdots x_{\xi_{m-1}} = \left\lvert x_{\xi_0} \right\rangle 
\otimes\!\cdots\otimes\!\left\lvert x_{\xi_{m-1}}\right\rangle; 
\,~ \xi=0,\ldots,n\!-\!1. 
\end{align}
$\,$\\[-.5in]
$\,$\\[-.5in]
\begin{theorem}[$\bm{\Phi}$ existence]\label{phiexits} 
The measure $\Phi$ exists$.$
\end{theorem}
$\,$\\[-.25in]
\noindent {\it Proof.}~~Let $\Phi\!:\mathcal{A}^{\otimes m}\!\!\longrightarrow\!\mathbb{R}_+$ satisfy the natural property for every observable ${\rm\bf \Gamma},$ resp$.$ ${\rm\bf \Gamma}.$ Then there is a unique state density ${\rm\bf T},$ resp$.$ ${\rm\bf T},$ such that$,$ just as in the continuous state$,$ 
\\[-.225in] 
\begin{align}
\Phi({\rm\bf\Gamma})\,= 
\sum_x{\rm\bf P}_{\!\!x}\,tr(\widetilde{\rm\bf T}_{\!x}\,{\rm\bf\Gamma}_{\!x}) 
\,=\sum_x tr({\rm\bf T}_{\!x}\,{\rm\bf\Gamma}_{\!x}) 
\,=\,tr({\rm\bf T}\,{\rm\bf\Gamma})
\,=\,\mathbb{E}^{gen}[{\rm\bf\Gamma}]
\hskip.2in\big|\hskip.1in 
\widetilde{\rm\bf T}_{\!x} = \ket{x}\!\bra{x} 
\nonumber 
\end{align}
\\[-.45in]
\begin{align}
{\rm\bf T}=\sum_x{\rm\bf T}_{\!x}\,
=\sum_x{\rm\bf P}_{\!\!x}\widetilde{\rm\bf T}_{\!x}
\hskip.125in\big|\hskip.075in 
{\rm\bf\Gamma}_{\!x}={\rm\bf T}^{-1/2}\,{\rm\bf T}_{\!x}\,{\rm\bf T}^{-1/2} 
=\,{\rm\bf T}^{-1/2}\,({\rm\bf P}_{\!x}\widetilde{\rm\bf T}_{\!x})\,{\rm\bf T}^{-1/2} 
\geqslant 0\,, 
\hskip.1in 
\sum_x{\rm\bf\Gamma}_{\!x}={\rm\bf I}_n.
\nonumber 
\end{align}
\\[-.175in] 
Hence$,$ the conclusion follows by Tab$.$~\ref{phi-table}$,$ under (\ref{phiexists})$.$ \hfill $\square$
$\,$\\[.1in]
{\bf Derivation.}~~If $n\!=\!2,$ i.e$.$ $\lvert\xi\rvert\!=\!3,$ $m\!=\!1,$ in particular$,$ for $\mathcal{A}^\otimes\!=\!(\ket{0}\!,\ket{1}),$ then $\Phi({\rm\bf\Gamma})\!=\!\mathbb{E}^{gen}[{\rm\bf\Gamma}]\!=\!\dfrac{5}{9}.$
\\[-.0in]
\begin{theorem}[optimality]\label{opttheorem} 
For an observable ${\rm\bf\Gamma}\!\mid\! \{{\rm\bf\Gamma}_{\!x}\}_{x\,\in\,\mathcal{X}},$ general $\mathbb{E}^{gen}[{\rm\bf\Gamma}],$ resp$.$ optimal $\mathbb{E}^{opt}[{\rm\bf\Gamma}],$ 
\\[-.2in]
\begin{align}
\mathbb{E}^{gen}[{\rm\bf\Gamma}] \,\geqslant\, 
\big(\mathbb{E}^{opt}[{\rm\bf\Gamma}]\big)^{\!2} 
\hskip.15in \text{i.e$.$}
\hskip.15in \mathbb{E}^{opt}[{\rm\bf\Gamma}] 
\leqslant\sqrt{\mathbb{E}^{gen}[{\rm\bf\Gamma}]} ~.
\nonumber 
\end{align}
\end{theorem}
$\,$\\[-.2in]
\noindent{\it Proof.}~~
\\[-.45in]
\begin{align}
\mathbb{E}^{gen}[{\rm\bf\Gamma}] 
\,=\,\sum_x tr\big( 
{\rm\bf T}_{\!x} 
\,{\rm\bf\Gamma}_{\!x}\big) 
\,=\,\sum_x tr\Big(\Big( 
\!{{\rm\bf T}}^{1/2} 
\,{\rm\bf\Gamma}_{\!x} 
\,{{\rm\bf T}}^{1/2} 
\Big)\Big( 
\!{{\rm\bf T}}^{-1/2} 
\,{\rm\bf T}_{\!x}
\,{{\rm\bf T}}^{-1/2} 
\Big)\Big) 
\nonumber 
\end{align}
\\[-.45in] 
\begin{align}
\leqslant\, 
\sum_x\,\Big(tr\Big( 
\!{{\rm\bf T}}^{1/2} 
\,{\rm\bf\Gamma}_{\!x} 
\,{{\rm\bf T}}^{1/2} 
\,{\rm\bf\Gamma}_{\!x} 
\Big)\Big)^{\!1/2}\,\Big(tr\Big( 
\!{{\rm\bf T}}^{-1/2} 
\,{\rm\bf T}_{\!x}
\,{{\rm\bf T}}^{-1/2} 
\,{\rm\bf T}_{\!x} 
\Big)\Big)^{\!1/2}  
\hskip.1in \text{\circled{$\leqslant$}} 
\nonumber 
\end{align}
$\,$\\[-.2in]
so that by the Cauchy-Bunyakovsky-Schwarz inequality for matrix-space Hilbert-Schmidt product
\\[-.225in]
\begin{align}
|tr(AB)| ~ \leqslant ~ 
\sqrt{tr(A A^\dagger)} ~ \sqrt{tr(B B^\dagger)}
\nonumber 
\end{align}
\\[-.25in]
therefore
\\[-.4in]
\begin{align}
&\text{\circled{$\leqslant$}} \hskip.1in 
\bigg(\!\sum_xtr\Big(\!{{\rm\bf T}}^{1/2}\,{\rm\bf\Gamma}_{\!x} 
\,{{\rm\bf T}}^{1/2}\,{\rm\bf\Gamma}_{\!x}\!\Big)
\!\bigg)^{\!\!1/2} 
\,\bigg(\!\sum_x tr\Big(\!{{\rm\bf T}}^{-1/2} 
\,{\rm\bf T}_{\!x}
\,{{\rm\bf T}}^{-1/2} 
\,{\rm\bf T}_{\!x}\!\Big)
\!\bigg)^{\!\!1/2} 
\hskip.05in \leqslant ~\sqrt{\mathbb{E}^{gen}[{\rm\bf\Gamma}]} 
\nonumber
\end{align}
\\[-.35in]
where
\\[-.25in]
\begin{align}
tr\Big(\!{{\rm\bf T}}^{1/2}\,{\rm\bf\Gamma}_{\!x} 
\,{{\rm\bf T}}^{1/2}\,{\rm\bf\Gamma}_{\!x}\!\Big) 
\leqslant\,tr\big({\rm\bf T}\,{\rm\bf\Gamma}_{\!x}\big), 
\hskip.15in 
tr\Big(\!{{\rm\bf T}}^{1/4}\,{\rm\bf\Gamma}_{\!x}
\,{{\rm\bf T}}^{1/4}\,{{\rm\bf T}}^{1/4}
\,{\rm\bf\Gamma}_{\!x}\,{{\rm\bf T}}^{1/4}\!\Big) 
=\, 
tr\Big(\!{{\rm\bf T}}^{1/2}\,{\rm\bf\Gamma}_{\!x}
\,{{\rm\bf T}}^{1/2}\,{\rm\bf\Gamma}_{\!x}\!\Big). 
\hskip.2in \square
\nonumber 
\end{align}
\\[-.45in] 
\begin{theorem}[quantum tripartite]\label{tritheorem} 
In $\varepsilon$-rescaling of $\mathbb{R}^3$ on finite $\mathcal{X}$ tree tangle $T_1\!\otimes_\varepsilon\!T_2,$ Fig$.$ \ref{tangle-rescale}$,$
for \cite{vj85,vj87,vj89,rt90,rt91} on the noted $K\!=\!G(X),$ uniform-density $\lvert\mathcal{A}\rvert$-universal hash functions family $G\!:$ 
\\[-.25in]
\begin{align}
G \ni g\!: \,\mathcal{A}^{\otimes \displaystyle m}\cong\mathcal{X} \longrightarrow\mathcal{A}^{\otimes \displaystyle k} \hskip.1in \mid 
\hskip.05in \mathcal{A}^{\otimes \displaystyle \zeta}\cong\,G 
\nonumber
\end{align}
\\[-.275in]
for $h_{\min}(X|Q)\geqslant h_+,$ there exists the metric
\\[-.4in]
\begin{align}
\lvert\mathcal{A}\rvert^{-1}\left\lVert 
{\rm\bf T}_{\!KGQ}\,-\,\lvert\mathcal{A}\rvert^{\displaystyle-k}
\,{\rm\bf I}_K\otimes{\rm\bf T}_{\!GQ} 
\right\rVert_{\mathcal{L}(L^1(\gamma))} 
\,\leqslant ~  \lvert\mathcal{A}\rvert^{^{\displaystyle\!\Big(\!\!-\frac{h_+\!-\!k}{2}\Big)}} = \varepsilon \hskip.1in \Big| ~ 0\!<\!\varepsilon\!\leqslant\!1. 
\end{align}
\end{theorem}
$\,$\\[-.4in]
$\,$\\[-.25in]
{\it Proof.}~~Write$:$ 
\\[-.4in]
\begin{align}
{\rm\bf T}_{\!KGQ} = \hskip-.1in 
\sum_{\kappa\,\in\,\mathcal{A}^{\otimes k}} 
\hskip-.075in\ket{\kappa}\!\bra{\kappa} \otimes {\rm\bf T}_{\!GQ,x} 
\hskip.15in \Big| \hskip.05in 
{\rm\bf T}_{\!GQ,\kappa} 
=\,\lvert\mathcal{A}\rvert^{-\zeta}\sum_{g\,\in\,G}\ket{x}\!\bra{x}  
\otimes\hskip-.175in 
\sum_{x:\,g(x)\,=\,\kappa} 
\hskip-.15in {\rm\bf T}_{\!Q,x} 
\end{align}
\\[-.45in]
\begin{align}
{\rm\bf T}_{\!XQ} \otimes {\rm\bf T}_{\!G} 
=\, \lvert\mathcal{A}\rvert^{\displaystyle -\zeta} 
\,\sum_{x,\,g}\,\ket{x}\!\bra{x}  
\otimes\ket{g}\!\bra{g} 
\otimes {\rm\bf T}_{\!Q,x} 
\hskip.15in \Big| \hskip.05in 
{\rm\bf T}_{\!G} =\, \lvert\mathcal{A}\rvert^{\displaystyle-\zeta} 
\,\sum_g\,\ket{g}\!\bra{g} 
\end{align}
\\[-.25in]
where 
\\[-.4in]
\begin{align}
{\rm\bf T}_{\!XQ} \,= \hskip-.0in \sum_{x\,\in\,\mathcal{X}} 
\hskip-.0in {\rm\bf P}_{\!\!x} 
\ket{x}\!\bra{x} \otimes\widetilde{{\rm\bf T}}_{\!Q,x} 
\,~ = \hskip-.0in \sum_{x\,\in\,\mathcal{X}} 
\hskip-.0in \ket{x}\!\bra{x} \otimes {\rm\bf T}_{\!Q,x}
\hskip.2in \big|\hskip.1in 
{\rm\bf T}_{\!Q,x} = {\rm\bf P}_{\!\!x}\,\widetilde{{\rm\bf T}}_{\!Q,x}. 
\end{align}
\\[-.15in]
For additional register (state)$,$
\\[-.25in]
\begin{align}
\ket{0}\!\bra{0}\otimes {\rm\bf T}_{\!XQ}\otimes {\rm\bf T}_{\!G} 
=\, \lvert\mathcal{A}\rvert^{\displaystyle -\zeta} 
\,\sum_{x,\,g}\,\ket{0}\!\bra{0}
\otimes\ket{x}\!\bra{x}  
\otimes\ket{g}\!\bra{g} 
\otimes {\rm\bf T}_{\!Q,x}. 
\end{align}
\\[-.25in]
Furthermore$,$
\\[-.375in]
\begin{align}
{\rm\bf T}_{\!KXGQ} 
=\, \lvert\mathcal{A}\rvert^{\displaystyle -\zeta} 
\,\sum_{x,\,g}\,\ket{g(x)}\!\bra{g(x)}
\otimes\ket{x}\!\bra{x}  
\otimes\ket{g}\!\bra{g} 
\otimes {\rm\bf T}_{\!Q,x} 
\end{align}
\\[-.4in]
\begin{align}
{\rm\bf T}_{\!KGQ} 
=\, \lvert\mathcal{A}\rvert^{\displaystyle -\zeta} 
\,\sum_\kappa\,\sum_g\sum_{x:\,g(x)\,=\,\kappa}\hskip-.125in\ket{\kappa}\!\bra{\kappa}  
\otimes\ket{g}\!\bra{g} 
\otimes {\rm\bf T}_{\!Q,x} 
\end{align}
\\[-.25in] 
where 
\\[-.35in]
\begin{align}
\sum_g\sum_x 
\,\ket{g(x)}\!\bra{g(x)}
\otimes\ket{x}\!\bra{x}  
\otimes\ket{g}\!\bra{g} 
\hskip.1in = \hskip.1in 
\sum_g\sum_\kappa 
\,\ket{\kappa}\!\bra{\kappa}
\otimes\hskip-.175in 
\sum_{x:\,g(x)\,=\,\kappa}\hskip-.125in 
\ket{x}\!\bra{x}  
\otimes\ket{g}\!\bra{g}
\end{align}
\\[-.375in]
\begin{align}
{\rm\bf T}_{\!KXGQ} 
\,=\, \lvert\mathcal{A}\rvert^{\displaystyle -\zeta}
\,\sum_\kappa\,\ket{\kappa}\!\bra{\kappa} 
\otimes\!\sum_g\,\ket{g}\!\bra{g} 
\otimes\hskip-.175in 
\sum_{x:\,g(x)\,=\,\kappa} 
\hskip-.125in\ket{x}\!\bra{x} 
\otimes {\rm\bf T}_{\!Q,x}
\end{align}
\\[-.375in]
\begin{align}\hskip-.0in
{\rm\bf T}_{\!GQ}\,&=\, tr_K({\rm\bf T}_{\!KGQ}) 
\,=\,\sum_\kappa {\rm\bf T}_{\!GQ,\kappa} 
\,=\,\lvert\mathcal{A}\rvert^{-\zeta}\sum_x {\rm\bf T}_{\!Q,x} 
\,=\,\lvert\mathcal{A}\rvert^{-\zeta}\sum_g\ket{g}\!\bra{g}\sum_x {\rm\bf T}_{\!Q,x} 
\nonumber  
\\[-.0in]
&=\lvert\mathcal{A}\rvert^{-\zeta}\sum_g\ket{g}\!\bra{g}\!\otimes\!{\rm\bf T}_{\!Q} 
= \lvert\mathcal{A}\rvert^{-\zeta} {\rm\bf I}_G\!\otimes\!{\rm\bf T}_{\!Q} 
\end{align}
\\[-.4in]
\begin{align}
{\rm\bf T}_{\!KGQ} = 
\sum_{\kappa,\,g}\,\ket{\kappa}\!\bra{\kappa}\otimes\ket{g}\!\bra{g}\otimes {\rm\bf T}_{\!Q,\kappa g} 
\hskip.15in \Big| \hskip.05in 
{\rm\bf T}_{\!Q,\kappa g} =\, 
\lvert\mathcal{A}\rvert^{\displaystyle -\zeta}\hskip-.15in\sum_{x:\,g(x)\,=\,\kappa} \hskip-.15in {\rm\bf T}_{\!Q,x}.
\end{align}
\\[-.15in] 
That is$,$ for collision probability$,$ on (Eva$,$ et al$.$) a priory known state ${\rm\bf T}_{\!Q},$ then the bound
\\[-.15in]
\begin{align}\hskip-.45in 
\lvert\mathcal{A}\rvert^{-1}\left\lVert{\rm\bf T}_{\!KGQ} 
\,-\,\lvert\mathcal{A}\rvert^{-k-\zeta} 
\,{\rm\bf I}_{KG}\otimes{\rm\bf T}_{\!Q}\right\rVert_{\mathcal{L}(L^1(\gamma))} 
=\,\lvert\mathcal{A}\rvert^{-1}\sum_{\kappa,\,g} 
\left\lVert {\rm\bf T}_{\!KG}(\kappa,g) 
\,-\,\lvert\mathcal{A}\rvert^{\displaystyle (-k\!-\!\zeta)} {\rm\bf T}_{\!Q}\right\rVert_{\mathcal{L}(L^1(\gamma))} \hskip.05in \text{\circled{$\leqslant$}}
\end{align}
\\[-.125in] 
where$,$ for matrices $\tau,\,\Lambda;$ $A\!=\!\tau^{-1/4}\,\Lambda\,\tau^{-1/4},\,B\!=\!\tau^{1/2}\,\mid\,tr(\tau)\!=\!1,$ $\tau\!=\!\tau^\dagger,$ the Cauchy-Bunyakovsky
\\[-.2in]
\begin{align}
|tr(AB)| \,\leqslant\, (tr(AA^\dagger))^{1/2} \,(tr(BB^\dagger))^{1/2}~, \hskip.2in 
|tr(\Lambda)| \,\leqslant \Big(tr\Big(\tau^{-1/4}\,\Lambda\,\tau^{-1/4}\Big)^{\!2}\,\Big)^{\!1/2} 
\nonumber 
\end{align}
\\[-.325in]
implies
\\[-.225in]
\begin{align}\hskip-.8in 
\text{\circled{$\leqslant$}} 
\hskip.15in \lvert\mathcal{A}\rvert^{-1}\sum_{\kappa,\,g}\!\bigg(
\!tr\Big(\!{{\rm\bf T}}^{-1/4}_{\!Q}
\Big(\!{\rm\bf T}_{\!Q,\kappa g} 
- \lvert\mathcal{A}\rvert^{\displaystyle -k\!-\!\zeta}\,{\rm\bf T}_{\!Q}
\Big){{\rm\bf T}}^{-1/4}_{\!Q}\Big)^{\!2}\bigg)^{\!\!1/2} 
\!\!\leqslant\,\lvert\mathcal{A}\rvert^{\displaystyle(k\!+\!\zeta)/2\!-\!1} 
\bigg(\!\sum_{\kappa,\,g} 
tr\Big(\!{{\rm\bf T}}^{-1/4}_{\!Q}\,{\rm\bf T}_{\!Q,\kappa g} 
- \lvert\mathcal{A}\rvert^{\displaystyle -k\!-\!\zeta} 
\,{{\rm\bf T}}^{-1/4}_{\!Q}\Big)^{\!2}\bigg)^{\!\!1/2} 
\nonumber 
\end{align}
\\[-.4in]
\begin{align} 
= ~ \lvert\mathcal{A}\rvert^{\displaystyle k/2\!-\!1} 
\Big(\lvert\mathcal{A}\rvert^{\displaystyle \zeta} 
\sum_{\kappa,\,g} 
tr\Big(\!{\rm\bf T}_{\!Q,\kappa g}{{\rm\bf T}}^{-1/4}_{\!Q} 
\,{\rm\bf T}_{\!Q,\kappa g}{{\rm\bf T}}^{-1/4}_{\!Q} 
\Big)\,-\,\lvert\mathcal{A}\rvert^{\displaystyle -k}\,\Big)^{\!\!1/2} \hskip.1in (\ast\ast)
\end{align}
\\[-.15in]
for $\sum_{\kappa,\,g}{\rm\bf T}_{\!Q,\kappa g}\!=1.$ That is$,$ (Eva$,$ et al$.$) collision probability ${\rm\bf P}_{\!col}$ is given by trace ``$tr$'' formula$.$ 
\\[.1in]
Precisely$,$ for the $K$ and $G$ joint collision probability ${\rm\bf P}_{\!col},$ with respect to expectation $\mathbb{E}^{gen}\!:$ 
\\[-.2in]
\begin{align}
\hskip-.2in {\rm\bf P}_{\!col} 
&=\lvert\mathcal{A}\rvert^{\displaystyle-2\zeta}
\,\sum_{\kappa,\,g}\hskip-.2in
\sum_{\substack{x,\,x':\\[.025in] g(x)\,=\,g(x')\,=\,\kappa}} 
\hskip-.275in 
tr\Big(\underbrace{\!\!{\rm\bf T}_{\!Q,x} 
{{\rm\bf T}}^{-1/2}_{\!Q} 
\,{\rm\bf T}_{\!Q,x} 
{{\rm\bf T}}^{-1/2}_{\!Q}\!}_{\displaystyle \#\#}\Big)
=
\lvert\mathcal{A}\rvert^{\displaystyle-2\zeta} 
\,\sum_{\kappa,\,g}\!\left(\hskip-.05in 
\begin{array}{l} 
\displaystyle 
\sum_{\substack{x\,\neq\,x':\\[.025in] 
g(x)\,=\,g(x')\,=\,\kappa}} 
\hskip-.275in \#\# \hskip.1in + 
\hskip-.15in \sum_{\substack{x\,=\,x':\\[.05in] 
g(x)\,=\,g(x')\,=\,\kappa}} 
\hskip-.275in \#\# 
\end{array}
\hskip-.075in\right) 
\\[.025in]
&=\,\lvert\mathcal{A}\rvert^{\displaystyle-(k\!+\!\zeta)}
\!\sum_{x\,\neq\,x'}
\!tr\Big(\!{\rm\bf T}_{\!Q,x}{{\rm\bf T}}^{-1/2}_{\!Q} 
\,{\rm\bf T}_{\!Q,x'}{{\rm\bf T}}^{-1/2}_{\!Q} 
\Big) ~ + ~ \lvert\mathcal{A}\rvert^{\displaystyle-\zeta}
\,\sum_x 
\,tr \Big(\!{\rm\bf T}_{\!Q,x} {{\rm\bf T}}^{-1/2}_{\!Q} 
\,{\rm\bf T}_{\!Q,x}{{\rm\bf T}}^{-1/2}_{\!Q} 
\Big) 
\\[-.025in]
&\leqslant\,\lvert\mathcal{A}\rvert^{\displaystyle-k\!-\!\zeta} 
\,+\,\lvert\mathcal{A}\rvert^{\displaystyle-\zeta\!-\!h_+} 
\hskip.15in \big| \hskip.1in 
{{\rm\bf T}}^{-1/2}_{\!Q}\,{\rm\bf T}_{\!Q,x'}
{{\rm\bf T}}^{-1/2}_{\!Q} = {\rm\bf\Gamma}_{x'}, 
\hskip.15in 
{{\rm\bf T}}^{-1/2}_{\!Q}\,{\rm\bf T}_{\!Q,x} 
{{\rm\bf T}}^{-1/2}_{\!Q} = {\rm\bf\Gamma}_{\!x}
\end{align} 
\\[-.2in]
that is$,$
\\[-.4in] 
\begin{align}
\sum_x 
\,tr\big({\rm\bf T}_{\!Q,x}
\,{\rm\bf\Gamma}_{\!x} 
\big) ~ =\, \mathbb{E}^{gen}[{\rm\bf\Gamma}] 
\,\leqslant\,\mathbb{E}^{opt}[{\rm\bf\Gamma}] 
\,=\,\lvert\mathcal{A}\rvert^{\displaystyle -h_{\min}(X|Q)} 
\,\leqslant\,\lvert\mathcal{A}\rvert^{-h_+}. \hskip.3in \square
\end{align}
\\[-.2in]
By $(\ast\ast),$ then
\\[-.25in]
\begin{align}
\lvert\mathcal{A}\rvert^{-1}\left\lVert \cdots \right\rVert_{\mathcal{L}(L^1(\gamma))}\,\leqslant\, 
\lvert\mathcal{A}\rvert^{\displaystyle k/2\!-\!1\!-\!k/2} 
\,\leqslant\, 
\lvert\mathcal{A}\rvert^{\displaystyle (k\!-\!h_+)/2-1} 
\cdot \lvert\mathcal{A}\rvert^{\displaystyle (k\!-\!h_+)/2}. \hskip.3in \square
\end{align}
\\[-.2in]
{\it Remark.}~~As a result$,$ the tree process hardness is proved for bb84 \cite{bb84} protocol strength$.$\vskip-.2in 

\subsection*{Addendum to asymptotic shape approximation with random trees} $\,$\vskip-.2in

\noindent Approximating locally finite tree shapes is a growing area of interest$;$ using the trace-metric process methods$,$ given in this article$,$ can provide rather faster convergence techniques for analytic  asymptotics$,$ on trees of same topology (that is, ``shape'')$,$ given by limit on increasing family of real (rooted$,$ finite) trees which appear in the construction of dynamical approximation techniques$,$ on interval $[0,1].$ Moreover$,$ the process discussed is just a random elements of the space $C[0,1]$ of continuous functions on $[0,1],$ by a finite field set with respect to supersymmetric (quantum) basis structure measures$,$ under suitable metric of stepwise intervals on right joint state asymptotics$;$ all other processes have duality to the given fundamental nature of the interpolation (construction)$.$ A necessary task is to show a {\it threshold of uniform convergence} for the equation (\ref{eq:interpolation})$.$ Such endeavor can be useful in tree asymptotics of braid (string) invariants for knot dualities$:$ Fig$.$ \ref{tangle-rescale}$.$
$\,$\\[-.1in]
\begin{figure}[H]
\centering ~ \hskip-0.in 
\includegraphics[width=.25\textwidth,angle=0,origin=c]
{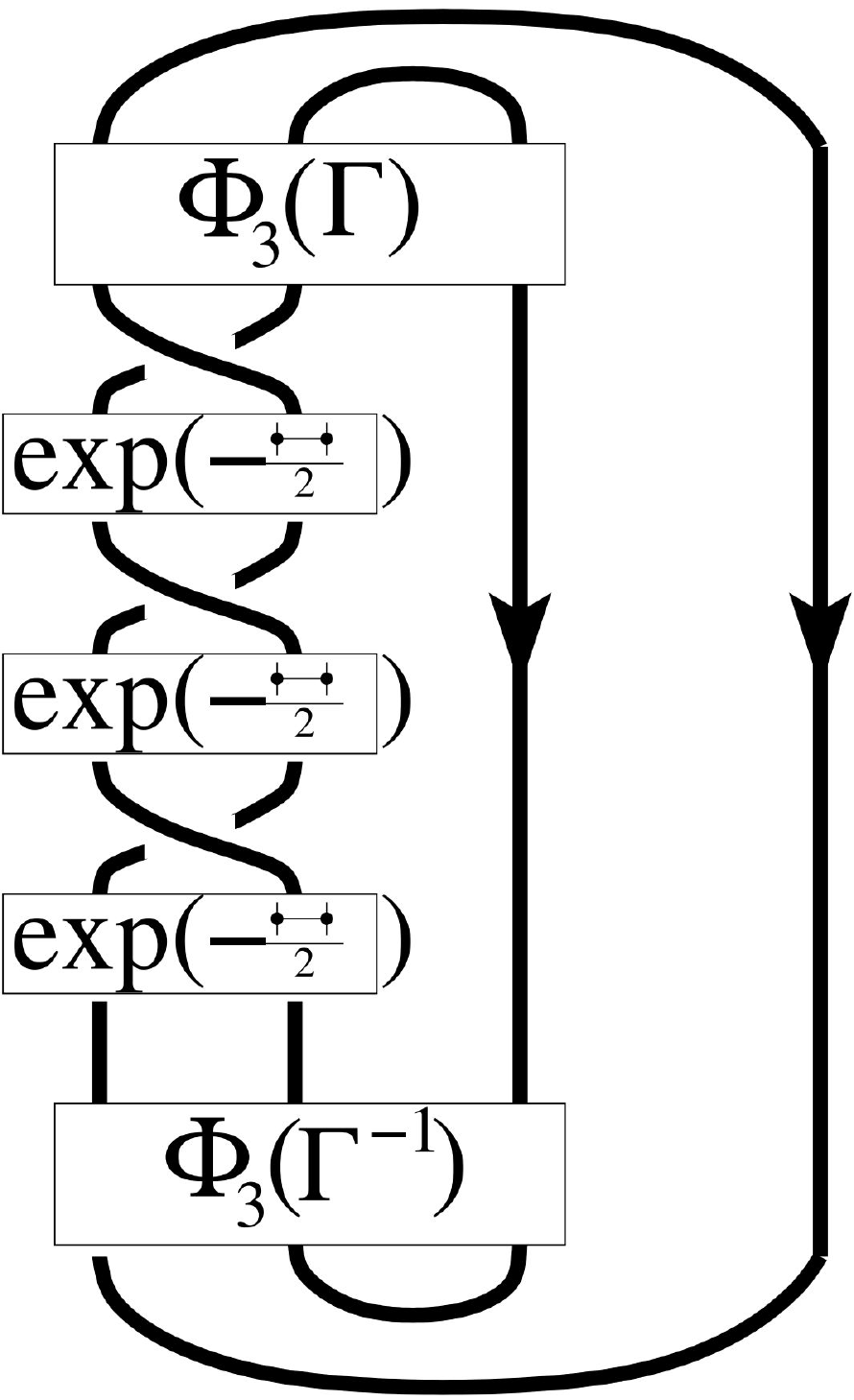}\hskip.3in 
\includegraphics[width=.4\textwidth,angle=0,origin=c]
{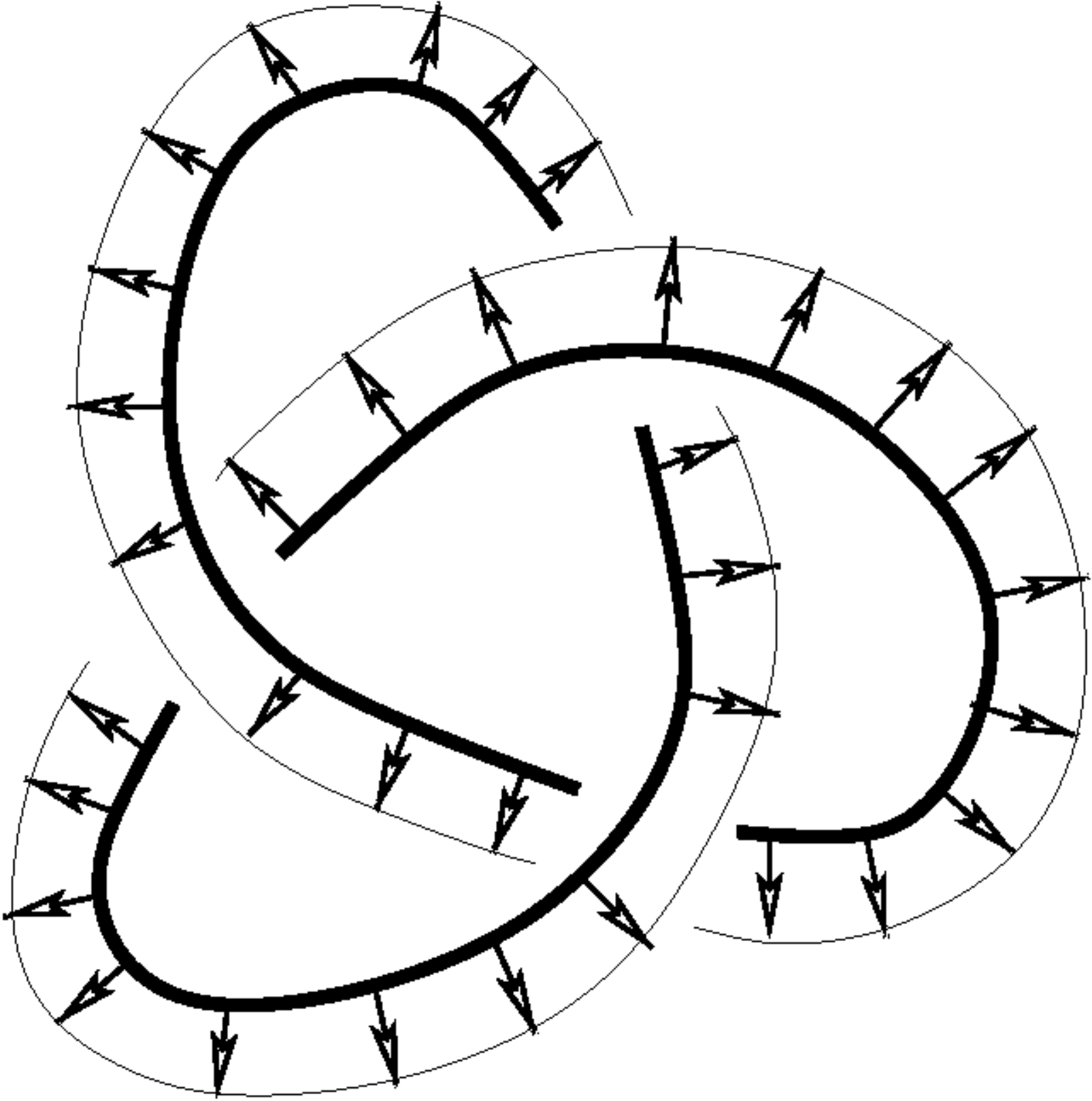}
\hskip-0in\,\\[-.05in]
\caption{Tree measure construction (left) for orientable-manifold$,$ framed-knot spin-structure duality (right)$.$}\label{string-algebra} \vskip-.075in 
\end{figure}
$\,$\\[-.45in]

\setcounter{section}{0}
\section*{Appendix}\label{appendix} \vskip-.05in 
\subsection{Closure of the $\bm{\lvert\mathcal{A}^{\otimes m}\rvert}$ normalization}\label{double}

The $1$ is trivial$.$ For $(2\!\times\!2),$ i.e$.$ $\lvert\mathcal{A}^{\otimes m}\rvert\!=\!2,$ write $t\Sigma$ and $t^{-1}\Sigma^{-1}$ (by cofactor or symmetric group)$:$
\\[-.175in]
\begin{align}
t\Sigma\,=\, t\left(\hskip-.05in
\begin{array}{cc}
\Sigma_{11} 
&\rho_{12}\sqrt{\Sigma_{11}\Sigma_{22}} 
\\[.075in]
\rho_{12}\sqrt{\Sigma_{11}\Sigma_{22}} 
&\Sigma_{22}
\end{array}
\hskip-.05in\right), 
\hskip.2in t^{-1}\Sigma^{-1}\,=\, \frac{1}{t\Sigma_{11}\Sigma_{22}\left(1-\rho_{12}^2\right)} 
\left(\hskip-.05in
\begin{array}{cc}
\Sigma_{22} 
&-\rho_{12}\sqrt{\Sigma_{11}\Sigma_{22}} 
\\[.075in]
-\rho_{12}\sqrt{\Sigma_{11}\Sigma_{22}} 
&\Sigma_{11}
\end{array}
\hskip-.05in\right) 
\nonumber 
\end{align}
\\[-.1in]
then we require the closure 
\\[-.2in]
\begin{align}\hskip-.0in 
\int_{\mathbb{R}}\int_{\mathbb{R}}\,\frac{1}{2\pi t\sqrt{\Sigma_{11}\Sigma_{22}\left(1-\rho_{12}^2\right)}} 
\exp\bigg\{\!\!-\frac{1}{2t}\bigg(
\frac{\left(x_1-\mu_1\right)^2}
{\left(1-\rho_{12}^2\right)\Sigma_{11}}\,-\,
2\frac{\left(x_1-\mu_1\right)\left(x_2-\mu_2\right)\rho_{12}}
{\left(1-\rho_{12}^2\right)\sqrt{\Sigma_{11}\Sigma_{22}}}\,+\,
\frac{\left(x_2-\mu_2\right)^2}
{\left(1-\rho_{12}^2\right)\Sigma_{22}}
\bigg)\bigg\}\,dx_1dx_2
\nonumber
\end{align}
\\[-.4in]
\begin{align}\hskip-.0in 
=\int_{\mathbb{R}}\int_{\mathbb{R}}\,\frac{1}{2\pi t\sqrt{\Sigma_{11}\Sigma_{22}\left(1-\rho_{12}^2\right)}} 
\exp\bigg\{\!\!-\frac{1}{2t}\bigg(
\bigg(\!\frac{1}{\sqrt{1-\rho_{12}^2}}\!\left(\!\frac{x_1-\mu_1}
{\sqrt{\Sigma_{11}}}\,-\,
\frac{\left(x_2-\mu_2\right)\rho_{12}}
{\sqrt{\Sigma_{22}}}\!\right)\!\!\bigg)^{\!\!2}+
\left(\!\frac{x_2-\mu_2}{\sqrt{\Sigma_{22}}}\!\right)^{\!\!2}
\,\bigg)\bigg\}\,dx_1dx_2
\nonumber
\end{align}
\\[-.1in] 
which is the closure of the origin-centered symmetry 
\\[-.25in] 
\begin{align}
\mathbb{S}^{n-1}\!=\!\bigg\{y_1,\ldots,y_n\in\mathbb{R} 
~ \bigg| ~ 
0<r=\sqrt{\sum^n_{i=1}y^2_i} ~ \bigg\}
\nonumber
\end{align}
\\[-.25in] 
for $n\!=\!2.$ 
\\[.1in] 
The closure of $2^m\!\mid\!m>1$ follows by induction from the $(3\!\times\!3)\!:$ 
\\[-.15in]
\begin{align}
t\Sigma\,=\, t \left(\hskip-.05in 
\begin{array}{ccc}
\Sigma_{11}
&\rho_{12}\sqrt{\Sigma_{11}\Sigma_{22}} 
&\rho_{13}\sqrt{\Sigma_{11}\Sigma_{33}} 
\\[.075in]
\rho_{12}\sqrt{\Sigma_{11}\Sigma_{22}}
&\Sigma_{22}
&\rho_{23}\sqrt{\Sigma_{22}\Sigma_{33}} 
\\[.075in]
\rho_{13}\sqrt{\Sigma_{11}\Sigma_{33}} 
&\rho_{23}\sqrt{\Sigma_{22}\Sigma_{33}} 
&\Sigma_{33}
\end{array}
\hskip-.05in\right),
\nonumber 
\end{align}
\\[-.3in]
\begin{align}
t^{-1}\Sigma^{-1}\,=\, &\frac{1}{t\Sigma_{11}\Sigma_{22}\Sigma_{33}\left(1-\rho_{12}^2-\rho_{13}^2-\rho_{23}^2+2\rho_{12}\rho_{13}\rho_{23}\right)} ~ \times 
\nonumber 
\\[.1in]
&\times\left(\hskip-.05in
\begin{array}{ccc}
\left(1-\rho_{23}^2\right)\Sigma_{22}\Sigma_{33} 
&\left(\rho_{13}\rho_{23}-\rho_{12}\right)
\Sigma_{33}\sqrt{\Sigma_{11}\Sigma_{22}} 
&\left(\rho_{12}\rho_{23}-\rho_{13}\right)\Sigma_{22} 
\sqrt{\Sigma_{11}\Sigma_{33}} 
\\[.075in]
\left(\rho_{13}\rho_{23}-\rho_{12}\right)
\Sigma_{33}\sqrt{\Sigma_{11}\Sigma_{22}} 
&\left(1-\rho_{13}^2\right)\Sigma_{11}\Sigma_{33} 
&\left(\rho_{12}\rho_{13}-\rho_{23}\right)\Sigma_{11} 
\sqrt{\Sigma_{22}\Sigma_{33}} 
\\[.075in]
\left(\rho_{12}\rho_{23}-\rho_{13}\right)\Sigma_{22}
\sqrt{\Sigma_{11}\Sigma_{33}} 
&\left(\rho_{12}\rho_{13}-\rho_{23}\right)\Sigma_{11} 
\sqrt{\Sigma_{22}\Sigma_{33}} 
&\left(1-\rho_{12}^2\right)\Sigma_{11}\Sigma_{22}
\end{array}
\hskip-.05in\right). 
\nonumber 
\end{align}
\\[-.05in]
In general$,$ by symmetric group $\mathcal{S}_n,$
\\[-.2in]
\begin{align} 
\det(\Sigma^{n\times n})\,=\, 
\sum_{\sigma\,\in\,\mathcal{S}_n}(-1)^{t(\sigma)}
\prod^n_{k=1}\Sigma_{k,\sigma(k)} 
\hskip.1in \bigg| ~ \Sigma_{i,j}=\Sigma_{ij}; 
\hskip.1in t(\sigma):= (\sigma(1),\ldots,\sigma(n))\longrightarrow(1,\ldots,n) 
\nonumber 
\end{align}
\\[-.125in]
for number of transpositions $t(\sigma)$ in $(\sigma(1),\ldots,\sigma(n))\longrightarrow(1,\ldots,n),$ for $\mathcal{S}_n$ automorphism $\sigma.$
\\[-.1in]
\begin{table}[H]
\begin{center}
\begin{tabular}{c|c}
\hline
\\[-.15in]
$n=2$ 
&$\Sigma_{11}\Sigma_{22}(1\!-\!\rho^2_{12})$
\\[.025in]
\hline
\\[-.15in]
\multirow{1}{*}{$n=3$} 
&$\Sigma_{11}\Sigma_{22}\Sigma_{33}(1\!-\!\rho^2_{12}\!-\!\rho^2_{13}\!-\!\rho^2_{23}\!+\!2\rho_{12}\rho_{13}\rho_{23})$
\\[.025in]
\hline
\\[-.15in]
\multirow{2}{*}{$n=4$} 
&$\Sigma_{11}\Sigma_{22}\Sigma_{33}\Sigma_{44}(1\!-\!\rho^2_{12}\!-\!\rho^2_{13}\!-\!\rho^2_{14}\!-\!\rho^2_{23}\!-\!\rho^2_{24}\!-\!\rho^2_{34}\!+\!\rho^2_{12}\rho^2_{34}\!+\!\rho^2_{13}\rho^2_{24}\!+\!\rho^2_{14}\rho^2_{23}$
\\[.025in]
&$+2\rho_{12}\rho_{13}\rho_{23}\!+\!2\rho_{12}\rho_{14}\rho_{24}\!+\!2\rho_{13}\rho_{14}\rho_{34}\!+\!2\rho_{23}\rho_{24}\rho_{34}\!-\!2\rho_{13}\rho_{14}\rho_{23}\rho_{24}\!-\!2\rho_{12}\rho_{14}\rho_{23}\rho_{24}\!-\!2\rho_{12}\rho_{13}\rho_{23}\rho_{34})$
\\[.025in]
\hline
\end{tabular}
\end{center}$\,$\vskip-.2in
\caption{Table of $\det(\Sigma)$}
\label{tab:multicol}
\end{table}
$\,$\\[-.2in]
\subsection{$\bm{\Phi}$ algebra}\label{phiexists}
$\,$\vskip-.25in

\begin{table}[H]
\hskip-.8in\begin{tabular}{|l|}
\hline
\\[-.125in] 
$\widetilde{\rm\bf T}_{\!\xi} =\, \ket{x_{\xi_0}\!\!\otimes\!\cdots\!\otimes\!x_{\xi_{m-1}}}\!\bra{x_{\xi_0}\!\!\otimes\!\cdots\!\otimes\!x_{\xi_{m-1}}},$ $\forall\,\xi\!=\!0,\ldots,{n\!-\!1},$ 
\,i.e$.$ projection operator into eigenvalue $\lambda_{\xi_0\cdots\xi_{m-1}}$ eigenspace
\\[.05in]
\hline
\\[-.125in]
$\widetilde{\rm\bf T}_{\!n} =\, \sum^{n-1}_{\xi=0}\widetilde{\rm\bf T}_{\!\xi} 
\,= \dfrac{1}{n}{{\rm\bf I}_n}$ 
\, i.e$.$ mixed state operator$;$ 
\hskip.2in 
${\rm\bf P}_{\!1} 
=\, \cdots 
=\, {\rm\bf P}_{\!n+1} 
=\, \dfrac{1}{n\!+\!1}$ 
\\[.1in]
\hline
\\[-.125in]
${\rm\bf T}=\sum^n_{\xi=0}{\rm\bf T}_{\!\xi} 
\,=\dfrac{1}{n\!+\!1}\Big(\!\sum^{n-1}_{\xi=0}\ket{x_{\xi_0}\!\!\otimes\!\cdots\!\otimes\!x_{\xi_{m-1}}}\!\bra{x_{\xi_0}\!\!\otimes\!\cdots\!\otimes\!x_{\xi_{m-1}}} 
\,+\,\dfrac{1}{n} \,{\rm\bf I}_n\!\Big) 
\,=\dfrac{1}{n\!+\!1}\Big(1\!+\!\dfrac{1}{n}\Big){\rm\bf I}_n 
~=\dfrac{1}{n} {\rm\bf I}_n;$ \hskip.1in 
${{\rm\bf T}}^{-1/2} =\, \sqrt{n} ~ {\rm\bf I}_n$
\\[.1in]
\hline
\\[-.125in]
${\rm\bf\Gamma}_{\!\!\xi} =\, n\!\cdot\!\dfrac{1}{n\!+\!1}\ket{x_{\xi_0}\!\!\otimes\!\cdots\!\otimes\!x_{\xi_{m-1}}}\!\bra{x_{\xi_0}\!\!\otimes\!\cdots\!\otimes\!x_{\xi_{m-1}}},$ ~$\forall\,\xi\!=\!0,\ldots,n\!-\!1;$ \hskip.1in 
${\rm\bf\Gamma}_{\!\!n} =\, 
n\!\cdot\!\dfrac{1}{n\!+\!1}\!\cdot\!\dfrac{1}{n} {\rm\bf I}_n 
\,=\, \dfrac{1}{n\!+\!1} {\rm\bf I}_n;$ \hskip.1in 
$\sum^n_{\xi=0} {\rm\bf\Gamma}_{\!\!\xi} 
\,=\, {\rm\bf I}_n$ 
\\[.1in]
\hline
\\[-.15in]
$\Phi({\rm\bf\Gamma})
=\,\dfrac{1}{n\!+\!1}\bigg(\!\!
\sum^{n-1}_{\xi=0}tr\Big(\!\ket{x_{\xi_0}\!\!\otimes\!\cdots\!\otimes\!x_{\xi_{m-1}}}\!\bra{x_{\xi_0}\!\!\otimes\!\cdots\!\otimes\!x_{\xi_{m-1}}}\!\cdot\! 
\dfrac{n}{n\!+\!1} \ket{x_{\xi_0}\!\!\otimes\!\cdots\!\otimes\!x_{\xi_{m-1}}}\!\bra{x_{\xi_0}\!\!\otimes\!\cdots\!\otimes\!x_{\xi_{m-1}}}\!\Big)+\,
tr\Big(\dfrac{1}{n} {\rm\bf I}_n\!\cdot\!
\dfrac{1}{n\!+\!1} {\rm\bf I}_n\!\Big)
\!\bigg)$ 
\\[.0in]
\\[-.15in]
$\hskip.2in =\,\dfrac{1}{n\!+\!1}\Big(\dfrac{n}{n\!+\!1}\sum^{n-1}_{\xi=0}\scaleone{2ex}_{\hskip-.0125in\xi}\big((\delta_{\xi\eta})|^n_{\xi,\eta=1}\big)\,+\,\dfrac{1}{n\!+\!1}\Big) 
\,=\,\dfrac{1}{n\!+\!1}\Big(n\dfrac{n}{n\!+\!1}+\dfrac{1}{n\!+\!1}\Big) 
\,=\,\dfrac{n^2\!+\!1}{(n\!+\!1)^2}\,=\,\mathbb{E}^{gen}[{\rm\bf\Gamma}]$
\\[.1in]
\hline
\end{tabular}\vskip-.05in
\caption{\hskip-.0in An existence of $\Phi$ for $\lvert\mathcal{A}^{\otimes m}\rvert\!+\!1$ quantum states$,$ i.e$.$ $(\lvert\mathcal{A}^{\otimes m}\rvert\!+\!1)$-mixture$,$ including mixed state$.$\hskip-.0in}\label{phi-table}
\end{table}
$\,$\\[-.3in]

\bibliographystyle{abbrvnat}
\nocite{*}
\bibliography{quantumcrypto-ref}

\section*{Acknowledgment} 
Research and authors were supported by the Lynn Bit Foundation in State of California$.$

\end{document}